\documentstyle[11pt,dissertation,rotate,sloppy]{USCthesis_mod}

\newfam\frakfam
\font\teneufm=eufm10     \textfont\frakfam=\teneufm
\font\seveneufm=eufm7    \scriptfont\frakfam=\seveneufm
\font\fiveeufm=eufm5     \scriptscriptfont\frakfam=\fiveeufm
\def\frak{\fam\frakfam \teneufm} 

\newfam\msbfam
\font\tenmsb=msbm10      \textfont\msbfam=\tenmsb
\font\sevenmsb=msbm7     \scriptfont\msbfam=\sevenmsb
\font\fivemsb=msbm5      \scriptscriptfont\msbfam=\fivemsb
\def\Bbb{\fam\msbfam \tenmsb}

\newcommand{\ie}{\mbox{\em i.e.}} 
\newcommand{\fq}{\mbox{${\Bbb F}_q\,$}}
\newcommand{\fp}{\mbox{${\Bbb F}_p\,$}}
\newcommand{\clo}{\mbox{${\overline {\Bbb F}}_p\,$}}
\newcommand{\nq}{\mbox{${\Bbb F}_{q^n}\,$}}
\newcommand{\sq}{\mbox{${\Bbb F}_{q^2}\,$}}
\newcommand{\sn}{\mbox{${\frak S}_n\,$}}
\newcommand{\op}{\mbox{${\mathcal O}_P\,$}}
\renewcommand{\oe}{\mbox{${\mathcal O}_E\,$}}
\newcommand{\oep}{\mbox{${\mathcal O}_{E,P}\,$}}
\newcommand{\of}{\mbox{${\mathcal O}_F\,$}}
\newcommand{\pf}{\mbox{${\Bbb P}(F)\,$}}
\newcommand{\pe}{\mbox{${\Bbb P}(E)\,$}}
\newcommand{\pti}{\mbox{${\Bbb P}(T_i)\,$}}
\newcommand{\jth}{\mbox{$j^{th}\;$}}
\newcommand{\ith}{\mbox{$i^{th}\;$}}
\newcommand{\disc}{\mbox{$\frak d$}}
\newcommand{\degdiff}{\mbox{\rm degDiff}}
\newcommand{\gal}{\mbox{\rm Gal}}
\newcommand{\diff}{\mbox{$\frak D$}}

\newcommand{\otone}{\mbox{${\mathcal O}_{T_1}\,$}}

\newtheorem{theorem}{Theorem}[section]
\newtheorem{lemma}[theorem]{Lemma}
\newtheorem{proposition}[theorem]{Proposition}
\newtheorem{corollary}[theorem]{Corollary}
\newtheorem{definition}[theorem]{Definition}

\newtheorem{example}[theorem]{\sc Example}

\renewcommand{\labelenumi}{\rm (\roman{enumi})}

\begin{document}
\bibliographystyle{plain}
\title{ ON SPLITTING PLACES OF DEGREE ONE IN EXTENSIONS OF ALGEBRAIC FUNCTION FIELDS, TOWERS OF FUNCTION FIELDS MEETING ASYMPTOTIC BOUNDS, AND BASIS CONSTRUCTIONS FOR ALGEBRAIC-GEOMETRIC CODES}
\author{Vinay Deolalikar}
\committee{
{\bf Oral Presentation} : & \\*[0.1in]
\hspace*{0.3in}	Date: April 16, 1999 (Friday)   & \\*
\hspace*{0.3in}	Time: 2:00 pm - 3:00 pm & \\*
\hspace*{0.3in}	Place: EEB-539   & \\*[0.2in]
{\bf Guidance Committee}:  & \\*[0.1in]
\hspace*{0.3in}    Late Prof. Dennis Estes & (Mathematics)\\*
\hspace*{0.3in}    Prof. P. Vijay Kumar	 & (Chairman)\\*
\hspace*{0.3in}	   Prof. Wayne Raskind & (Mathematics)\\*
\hspace*{0.3in}	   Prof. Solomon Golomb & (EE-Systems)\\*
\hspace*{0.3in}    Prof. Robert Guralnick & (Mathematics)\\*
\hspace*{0.3in}	   Prof. Lloyd Welch & (EE-Systems)
}
\submitdate{May 1999}
\majorfield{Electrical Engineering}
\begin{preface}

\prefacesection{}

\begin{center}

\it

To the late Prof. Dennis Estes, my greatest teacher,

\bigskip
\bigskip

and

\bigskip
\bigskip

To my late parents 

My mother, Smt. Usha Deolalikar, who taught me the value of courage, compassion and simplicity.

My father,  Shri. Shrinivas Deolalikar, who taught me the worth of an honest day's toil, and to stand up for what I believe.

\end{center}

\prefacesection{Acknowledgements}

Words fail me as I attempt to describe my debt to the one person who, above everyone else, helped me shape this work. My greatest teacher, the late Prof. Dennis Estes. While many may have said it before me, I mean it fully when I say that without his constant help and encouragement, this work would not have seen the light of day. If I have at all done anything in this work that deserves mention, it is entirely due to him. 

When I came to USC as a graduate student in the Fall of 1995, all I knew was that I wished to learn mathematics. I had, a few months back, been introduced to the idea of a group. From there on, Prof. Estes became my teacher in a long journey, where at times he had more faith in me than I. He taught me most of the mathematics that I learnt in these last four years. His encouragement and total support during these years gave me something to lean upon at all times. And lean I did. At the very first hint of a problem - mathematical or otherwise - I went to him, and he never let me down. 

Very tragically, he departed on the first of February, 1999, after having done  so much for me that a few words could never do justice to it.  Suffices to say that he occupies a holy space in my heart, along with my parents.  

I want to profusely thank my advisor, Prof. Vijay Kumar, whose generosity made my stay at USC a pleasure.  I am grateful to him for so many things that I learnt from him. Also, he gave me complete freedom to chart the course of my graduate study. He let me make all the crucial decisions, while always being there to offer advice. I think that without such freedom, I would not have been able to motivate myself. I feel very fortunate that he gave me the opportunity to work with him.  I gained a lot from our association and I hope it continues well into the future. 

I am truly indebted to Joe Wetherell. He spent so much of his time trying to help me, especially during those testing times following Prof. Estes' demise. He painstakingly read my work, pointed out several flaws and how to bridge them, and helped immensely in putting this work in its final shape. Not only did he do this while he was here, but continued while on travel, sending me crucial corrections and observations by email. That he did all this when he really need not have, shows the strength of his character. I am fortunate to have him as my friend, and as a teacher. I hope to learn much from him in the days to come. 

I wish to add a special acknowledgement to the late Prof. Subhashish Nag, who passed away last year, and from whom I learnt much. When he was visiting USC in the summer of 1997, we would meet almost everyday and he would discuss mathematics, from which I benefitted greatly. 

I thank my committee - Prof. Solomon Golomb, Prof. Wayne Raskind, Prof. Robert Guralnick, and Prof. Lloyd Welch - for their all their help and cooperation.

Now come my friends. PC first, because he more than anyone else has always had faith in me, and has tried to inspire me when he has felt that I needed it. His selfless friendship is my great strength. Lisa, who is always there for me, through thick and thin. Her friendship has meant so much to me. Ranvir, who has been my friend for so many years, in which we saw many ups and downs. The great times we had when he visited last summer are etched in my memory.  Mahesh, for his unconditional support at all times. Gaurav, for his friendship. My best friends in India - Arvind, Vikas and Kaushik. Ranga, for playing a crucial role in my thesis writing and for his company during that lonely endeavour.

I thank my family - Akka,  Manik Maushi, Pradeep mama, Manju mami, Suneel mama and  Sharayu mami. They live far from here, but are in my heart. They have given me a place I can call my home. My sisters Shilpa, Swati and Shweta, for whom I reserve my fondest hopes and dreams. Last, but certainly not the least, I thank Appa for his unflinching support ever since I came to America for graduate study.

\begin{singlespace}
\tableofcontents
 \end{singlespace}

\prefacesection{Abstract}

In this work, we use the notion of ``symmetry'' of functions for an extension $K/L$ of finite fields to produce extensions of a function field $F/K$ in which almost all places of degree one split completely. Then we introduce the notion of ``quasi-symmetry'' of functions for $K/L$, and demonstrate its use in producing extensions of $F/K$ in which all places of degree one split completely. Using these techniques, we are able to restrict the ramification either to one chosen rational place, or entirely to non-rational places.

We then apply these methods to the related problem of building asymptotically good towers of function fields. We construct examples of towers of function fields in which all rational places split completely throughout the tower. We construct Abelian towers with this property also. We also generalize two existing examples of towers of function fields meeting the Drinfeld-Vladut bound, resulting in two infinite families of such towers.

We obtain results on the existence of Abelian extensions of arbitrary characteristic power degree in which all rational places split completely, as well as non-Abelian extensions of arbitrarily high degree in which the same holds true.

Furthermore, all of the above are done explicitly, \ie, we give generators for the extensions, and equations that they satisfy. 

We also construct an integral basis for a set of places in a tower of
function fields meeting the Drinfeld-Vladut bound using the discriminant
of the tower localized at each place. Thus we are able to obtain a basis
for a collection of functions that contains the set of regular functions
in this tower.  Regular functions are of interest in the theory of
error-correcting codes  as they lead to an explicit description of the
code associated to the tower by providing the code's generator matrix.

\prefacesection{Notation}

For symmetric polynomials:
\begin{tabbing}
\sn \hspace{1.6cm} \= the symmetric group on $n$ characters \\
${\bf s}_{n,i}(X)$ \> the \ith elementary symmetric polynomial on $n$ variables \\ 
$q$  \>  a power of a prime $p$ \\
${\Bbb F}_l$  \> the finite field of cardinality $l$ \\
$s_{n,i}(t)$ \> the \ith $(n,q)$-elementary symmetric polynomial 
\end{tabbing}

For function fields and their symmetric subfields:
\begin{tabbing}
$K$ \hspace{1.6cm} \= the finite field of cardinality $q^n$, where $n>1$ \\
$F/K$ \> an algebraic function field in one variable whose full field of constants is $K$ \\
$F_s$ \> the subfield of $F$ comprising $(n,q)$-symmetric functions \\
$F_s^\phi$ \> the subfield of $F_s$ comprising functions whose coefficients are from \fq \\
$F_{qs}$ \> the subfield of $F$ comprising $(n,q)$-quasi-symmetric functions \\
$F^{\phi}_{qs}$ \> the subfield of $F_{qs}$ comprising functions whose coefficients are from \fq \\
$U_{qs}$ \> the \nq-vector space of $(n,q)$-quasi-symmetric functions from \nq to \nq\\
$V_{qs}$ \> the \fq-vector space of $(n,q)$-quasi-symmetric functions from \nq to \nq with a \\
	\> polynomial representation with \fq coefficients\\
$E$ \> a finite separable extension of $F$, $E=F(y)$ where $\varphi(y) = 0$ 
for some irreducible  \\
 \> polynomial $\varphi[T] \in F[T]$ 
\end{tabbing}

For a generic function field $F$:
\begin{tabbing}
${\Bbb P}(F)$  \hspace{1.3cm} \=  the set of places of $F$ \\
$N(F)$  \> the number of places of degree one in $F$ \\
$N_m(F)$\> the number of places of degree $m, m>1$, in $F$ \\
$g(K)$ \> the genus of $F$   \\
$P$ \> a generic place in $F$  \\
$v_P$ \> the normalized discrete valuation associated with the place $P$ \\
\op \> the valuation ring of the place $P$ \\
$P'$ \> a generic place lying above $P$ in a finite separable extension of $F$\\ 
$e(P'|P)$ \> the ramification index for $P'$ over $P$ \\
$d(P'|P)$ \> the different exponent for $P'$ over $P$ 
\end{tabbing}

For the rational function field $K(x)$:
\begin{tabbing}
$P_\alpha$ \hspace{1.6cm} \= the place in $K(x)$ that is the unique zero of $x-\alpha, \; \alpha \in K$ \\
$P_\infty$ \> the place in $K(x)$ that is the unique pole of $x$ 
\end{tabbing}

For towers of function fields:
\begin{tabbing}
${\mathcal F}$ \hspace{1.8cm} \= a tower of function fields $F_i \subseteq F_2 \subseteq F_3 \ldots$ \\
$\lambda(\mathcal F)$ \> $\displaystyle\lim_{i \rightarrow \infty}(N(F_i)/g(F_i))$ 
\end{tabbing}

For a fixed extension $E/K$ of the function field $F/K$:
\begin{tabbing}
$N_{E/F}(.)$ \hspace{0.6cm} \= the norm of an element or divisor, from $E$ to $F$ \\
\oep \> the integral closure of \op in $E$, where $P \in \pf$ \\
$\diff(E/F)$ \> the different of $E/F$ given by $\diff(E/F) = \displaystyle\prod_{P' \in \pe}P'^{d(P'|P)}$, it is a divisor in $E$ \\
$\disc(E/F)$ \> the discriminant of $E/F$, it is a divisor in $F$ given by $N_{E/F}(\diff(E/F))$ \\
$\disc(S)$ \> the discriminant of the set $S$ of elements from $E$ 
\end{tabbing}

\prefacesection{Preface}

An algebraic function field $F/K$ in one variable with a finite field of
constants $K$ is a finite algebraic extension of $K(x)$, where $x \in F$
is transcendental over $K$, and $K$ is assumed algebraically closed in
$F$. We will also refer to these as just function fields. Let ${\mathcal
C}$ be a nonsingular projective algebraic curve defined over $K$. Let $E$
be the field of rational functions on ${\mathcal C}$ with coefficients
from $K$. Then $E$ is an algebraic function field. Let $N(E)$ and $g(E)$
denote the number of places of degree one (or {\em rational} places), and
the genus, respectively, of $E$. This is the same as the number of
rational points on ${\mathcal C}$, and the genus of ${\mathcal C}$,
respectively.
Indeed, the geometric study of a curve is equivalent to the algebraic study of its associated function field. Also, there seem to be certain advantages to the latter approach when the underlying field is not algebraically closed and, in particular, if it is finite. Throughout this report, we will find it convenient to talk in the purely algebraic language of function fields. It will be understood that the affine equation of the curve being discussed is that satisfied by the generator of the function field, viewed as an extension of the rational function field. 

Thus stated, the objects of our study are algebraic function fields with many rational places. In recent years, there has been a tremendous amount of attention paid to this broad area. The impetus for this interest came, interestingly, from a fundamental advance in coding theory. Around 1980, Goppa \cite{Gop1} discovered that function fields with many rational places over a finite field could be used to construct linear codes over that field. Further, the code parameters, such as dimension and minimum distance, could then be estimated using the Riemann-Roch theorem. The rate of information transmission using such codes, which came to be  known as Algebraic-Geometric codes, depended upon the ratio $N/g$ of the number of rational places in the function field to its genus. $N$ is also the length of the code. An excellent reference for Algebraic-Geometric codes, as well as for function fields in general, is \cite{Sti1}.

But function fields with many rational places have always been objects of intrinsic mathematical interest. There are, for instance, connections of such function fields to several deep results in arithmetical algebraic geometry, going back to the celebrated ``Weil conjectures.'' Indeed, the Hasse-Weil bound, which is the fundamental bound on the number of rational places in a function field in terms of its genus, follows immediately from the analogue of the Riemann hypothesis for function fields. For a history of the Weil conjectures, up to their final resolution in the 1970s, pl. refer to \cite{FreKie1}.

One approach to obtain function fields with many rational places is to split many rational places  completely in an extension of some function field, usually taken to be the rational function field itself. This is the main idea behind the class field theoretic approaches to this problem,  introduced by Serre \cite{Ser1,Ser2,Ser3}. However, it is not always easy to describe such extensions, and there is a fair bit of computation involved before one can completely describe such a class field. Furthermore, these constructions are not explicit in that generators and equations satisfied by them are not provided. For practical applications of such function fields, however, it is necessary that the constructions be explicit. To overcome this problem, some authors have used explicit class field theory - specifically Drinfeld modules - to construct function fields with many rational places. This technique, however, presupposes a fair bit of advanced mathematics, and is inaccessible to the beginner. Also, in this case, it is perhaps even harder to describe the resulting class fields.

Our aim then, was two-fold. Firstly, to construct a mathematical machinery which could be used to  construct infinite families of examples of extensions of function fields with almost all rational places splitting completely, thus yielding a high number of rational places in the extension. Moreover, we wanted to do this in a simple way that would make the solution look ``natural'' in some sense. Secondly, to be able to explain many  existing examples of function fields with many rational places in terms of a broader theory.

 To these ends, we were able to build a general theory to construct function fields over finite fields using the notions of ``symmetry'' and ``quasi-symmetry.'' Not only did it yield infinite families of extensions with complete splitting of all rational places, but it also provided us, along the way,  with a ``generalization'' of the Hermitian function field over \sq to odd degree extensions  of ${\Bbb F}_q$. We showed that the Hermitian function field is a very special case of a large class of function fields, which can be constructed using symmetric functions.

We then used these techniques to build towers of function fields with many rational places for their genus, and we succeeded in constructing examples in which all rational places split completely at each step of the tower. We were also able to construct Abelian towers with this property.

 Furthermore, all our methods yield explicit generators for each extension and equations satisfied by them. Moreover, there is absolutely no computation involved in providing these. In that sense, these objects emerge canonically from our methods.

To actually construct codes on any tower of function fields that meets the Drinfeld-Vladut bound, it is necessary to have a basis for the space of those functions at every stage of the tower which have poles only at a certain fixed divisor and nowhere else. Usually, this fixed divisor is taken to be a power of the place at infinity. Functions whose only poles are at infinity are said to be regular. The problem of finding a basis for the space of regular functions in the tower meeting the Drinfeld-Vladut bound presented in \cite{GarSti2}, remains an open one, and has invited intense research activity. But recently, we have arrived at partial results pointing to a solution to this problem. These results are described in this report. We are optimistic that a final solution to this problem will be found soon.
\begin{flushright}
Vinay Deolalikar \\
Los Angeles,  April 22, 1999.
\end{flushright}

\end{preface}



\chapter{Symmetry} \label{chapter:symmetry}

\section{Introduction} 

In this chapter, we study the problem of splitting ``almost all'' rational places in extensions of function fields. We describe some families of such extensions in which all except one rational place split completely. The other technique existing in the literature that can be used to split many rational places uses class field theory, and was introduced by Serre \cite{Ser1,Ser2,Ser3}. For practical applications of such function fields, however, it is necessary that the constructions be explicit, in that generators and equations satisfied by them should be provided. Our approach succeeds in providing explicit descriptions of new families of extensions of function fields in which almost all rational places split completely. 

In a sense, such constructions aim for a high $N(E)/[E:F]$ ratio, rather than $N(E)/g(E)$. However, in many cases, including a large number of known cases, such constructions also yield function fields with many rational places for their genus. In some of these, the function field attains known bounds (Hasse-Weil, Oesterle) on such behaviour as well. Further, the function field towers of \cite{GarSti1} and \cite{GarSti2} also iterate special cases of these constructions.

Another point of considerable importance is that a disproportionately large number of the known function fields with many rational places for their genus are concentrated at lower values of genera, and have small size of the underlying field of constants. In other words, it is easier to produce ``good'' function fields for low values of genera and over a small field of constants. It gets harder to do the same as the genus and the size of the underlying field increase. However, for applications to coding theory, one would like to have a large number of rational places in order to build long codes. The function field towers of  \cite{GarSti1} and \cite{GarSti2}, give us a structure to build long codes, but the problem of finding a basis for the vector spaces of regular functions on them is yet to be solved. The constructions provided in this report can provide examples of extensions of function fields $E/F$ where all rational places except one split completely, for arbitrarily high degree of extension $[E:F]$, and therefore, for arbitrarily high values of $N(E)$. But it must be borne in mind that the size of the underlying field of constants would also need to be suitably large.

Many of the existing  explicit constructions of function fields with many rational places work for a field of constants of square cardinality only.  The families described in this report, however,  can be constructed for all finite fields of cardinality not equal to their prime characteristic. Moreover, the methods used to construct them shed light on why the square cardinality constraint often arises in existing examples.

Finally, we provide a ``generalization'' of the Hermitian function field over $K$ where $K$ is an odd degree extension  of ${\Bbb F}_q$. For a long time, it has been known that the Hermitian function field has many unique properties, such as its maximality in the Hasse-Weil sense, and its large automorphism group. However there has hitherto not been a satisfactory generalization of this function field over fields of nonsquare cardinality. We seek to address this problem.

\section{Symmetric functions} \label{section:symmetric-functions}

Let $R$ be an integral domain and $\overline{R}$ its field of fractions. Consider the polynomial ring  in $n$ variables over $R$, given by  $R\,[X]= R\,[x_1,x_2,\ldots,x_n]$. The symmetric group \sn acts in a natural way on this ring by permuting the variables. 

\begin{definition} 
A polynomial ${\bf f}(X) \in R\,[X]$ is said to be symmetric if it is fixed under the action of \sn. If \sn is allowed to act on $\overline{R}(X)$ in the natural way, its fixed points will be called symmetric rational functions, or simply, symmetric functions. These form a subfield $\overline{R}(X)_s$  of $\overline{R}(x)$.
\end{definition}

Symmetric functions form a very elegant branch in the study of polynomials with several indeterminates. We recall here one of the classical results on symmetric functions, often known as {\em The fundamental theorem on symmetric functions}\footnote{A sharper result, called {\em The fundamental theorem on symmetric polynomials}, says that every symmetric polynomial can be written as a polynomial in the elementary symmetric functions \cite{BeaPie1}.} \cite{BeaPie1}.
\begin{theorem}[The fundamental theorem on symmetric functions]
$$\overline{R}(X)_s  = \overline{R}({\bf s}_{n,1}(X), {\bf s}_{n,2}(X), \ldots, {\bf s}_{n,n}(X)) $$
where $\{{\bf s}_{n,i}(X)\}_{1\leq i \leq n}$ are given below:
\begin{eqnarray*}
  {\bf s}_{n,1}(X) &=& \sum_{i=1}^{n} x_i, \\
  {\bf s}_{n,2}(X) &=& \sum_{i<j \atop 1 \leq i,j \leq n} x_ix_j, \\
     \vdots  & &  \vdots \\
  {\bf s}_{n,n}(X) &=& x_1x_2\ldots x_n.
\end{eqnarray*}
 ${\bf s}_{n,i}(X)$ is called the \ith elementary symmetric polynomial in $n$ variables.
\end{theorem}

We would like to apply these notions to the setting of polynomial rings in one variable over finite fields. We start with the following elementary proposition.
\begin{proposition} \label{proposition:extension-finite-fields}
Consider the extension of finite fields given by ${\Bbb F}_{q^n}/{\Bbb 
F}_q$. Then, the following hold:
\begin{enumerate}
\item This is a Galois extension of degree $n$. The Galois group of this extension, $G=\gal({\Bbb F}_{q^n}/{\Bbb F}_q)$ is cyclic, and is generated by $\phi \in G, \; \phi:\alpha \rightarrow \alpha^q$. This generating element is also called the ``frobenius'' of this extension.
\item The fixed field of the subgroup of $G$ generated by $\phi^k$ is ${\Bbb F}_{q^{\gcd(k,n)}}$.
\end{enumerate}     
\end{proposition}

\begin{definition}
For the extension ${\Bbb F}_{q^n}/{\Bbb F}_q$, we will evaluate the elementary symmetric polynomials (resp. symmetric functions) in ${\Bbb F}_{q^n}(X)$ at $(X)=(t,\phi(t),\ldots,\phi^{n-1}(t))=(t,t^q,\ldots,t^{q^{n-1}})$. These will be called the $(n,q)$-elementary symmetric polynomials (resp. $(n,q)$-symmetric functions). For ${\bf f}(X) \in \nq\!(X)$, we will denote ${\bf f}(t,t^q,\ldots,t^{q^{n-1}})$ by $f(t)$, or, when the context is clear, by $f$.
\end{definition} 

Thus the  $(n,q)$-elementary symmetric polynomials are the following:
\begin{eqnarray*}
s_{n,1}(t) &=& \sum_{0 \leq i \leq n-1} t^{q^i}, \\
s_{n,2}(t) &=& \sum_{i<j \atop 0 \leq i,j \leq n-1}t^{q^i}t^{q^j}, \\
    \vdots & & \vdots  \\
s_{n,n}(t) &=& t^{1+q+q^2+\ldots+q^{n-1}}.
\end{eqnarray*}

\begin{lemma} \label{lemma:fsphi-in-fq}
Let $f(t)$ be an $(n,q)$-symmetric function with coefficients from \fq and $\gamma \in \nq$. Then, we have that $f(\gamma) \in {\Bbb F}_q \cup \infty$.
\end{lemma}
{\bf Proof}. $f(t)$ can be written as a rational function in $\{s_{n,i}(t)\}_{1 \leq i \leq n}$, with coefficients from \fq.  By construction, each $(n,q)$-elementary symmetric polynomial is invariant under the operation of raising to the $q^{th}$ power, modulo $(t^{q^n} - t)$. In other words, each $(n,q)$-elementary symmetric polynomial, when restricted to \nq, is invariant under this operation. So are the coefficients, since they are chosen from \fq. The inclusion of infinity in the range comes since $\gamma$ may be a pole of $f$. \hfill $\Box$
 
There are $n$ $(n,q)$-elementary symmetric polynomials. $s_{n,1}$ and $s_{n,n}$  are also known as ``trace'' and ``norm,'' respectively. In the case of $n=2$, these are the only $(n,q)$-elementary symmetric polynomials. For $n \geq 3$, we have a greater choice, which has not been exploited in existing techniques to construct function fields with many rational places. The main thesis of this chapter is that these additional polynomials are indeed the more useful ones, as $n$ increases. In particular, we show that a ``generalization'' of the Hermitian function field can be obtained for $n\geq3$ by using $s_{n,2}$.

We end this section with some lemmas that will be used later. All of the following are valid for $1 \leq i \leq n$. 

\begin{lemma} \label{lemma:Sni-Sn(n-i)}
$$ \frac{s_{n,i}(t)}{s_{n,n}(t)} = s_{n,n-i}\left(\frac{1}{t}\right).$$
Thus, there is a bijection between the set of nonzero roots of $s_{n,i}(t)$ and that of $s_{n,n-i}(t)$. More precisely, $s_{n,i}(\alpha) = 0 \Leftrightarrow s_{n,n-i}(1/\alpha) = 0$, where $\alpha \neq 0$.
\end{lemma}
{\bf Proof}. Follows immediately from the relation between $s_{n,i}(t)$ and $s_{n,n-i}(t)$ and observing that the only root of $s_{n,n}(t)$ is zero itself. 

\begin{lemma} \label{lemma:derivative-of-Sni}
$$[s_{n,i}(t)]' = {[s_{n-1,i-1}(t)]}^q. $$
\end{lemma}
{\bf Proof}. The only terms that will contribute to the derivative are those whose exponent is coprime to $q$, \ie, in which $1$ is a summand. There are ${n-1 \choose i-1}$ such terms and they are all distinct.  Moreover, each such term, divided by $t$, is just the $q^{th}$ power of a term in $s_{n-1,i-1}(t)$. It is easy to see that this correspondence is bijective by noting that there are exactly ${n-1 \choose i-1}$ terms in $s_{n-1,i-1}(t)$ as well. \hfill $\Box$

\begin{lemma} \label{lemma:Sni-splits-partially}
The roots of $s_{n,i}(t)$ of multiplicity coprime to $p$ lie in \nq.
\end{lemma}
{\bf Proof}. The proof rests on the observation that 
\begin{eqnarray*}
 {s_{n,i}(t)}^q - s_{n,i}(t) &=& [t^{q^n} - t] {[s_{n-1,i-1}(t)]}^q, \\
			     &=& [t^{q^n} - t] {s_{n,i}(t)}', \\
\mbox{or,                  } s_{n,i}(t)[{s_{n,i}(t)}^{q-1} - 1] &=& [t^{q^n} - t] {s_{n,i}(t)}'. 
\end{eqnarray*}
Now consider a root $\alpha$ of multiplicity $m$ of $s_{n,i}(t)$ in the algebraic closure of ${\Bbb F}_q$, where $\gcd(m,p)=1$. This is not a root of the factor $[{s_{n,i}(t)}^{q-1} - 1]$ and it can only occur to a multiplicity $m-1$ in ${s_{n,i}(t)}'$. So it must appear to multiplicity $1$ in the factor $[t^{q^n} - t]$ on the RHS. Thus it must lie in ${\Bbb F}_{q^n}$.  \hfill $\Box$

\begin{lemma} 
Every root of $s_{n,i}(t)$ lies in a field ${\Bbb F}_{q^k}$, where $n-i+1 \leq k \leq n$.
\end{lemma}
{\bf Proof}. From Lemma~\ref{lemma:Sni-splits-partially}, we know that if a root of $s_{n,i}(t)$ is not also a root of $s_{n-1,i-1}(t)$, it must lie in \nq. If, on the other hand, it is a root of $s_{n-1,i-1}(t)$,  it either lies in ${\Bbb F}_{q^{n-1}}$ or is a root of $s_{n-2,i-2}$. Now we descend this way till we reach $s_{n-i+1,1}(t)$, all of whose roots lie in ${\Bbb F}_{q^{n-i+1}}$. Notice that $s_{n-i+1,1}(t)$ can have no multiple roots since its derivative is equal to $1$. \hfill $\Box$

\begin{lemma} 
Let $\alpha$ be a root of $s_{n,i}(t)$ of multiplicity $m>1.$ Further let
$\gcd(m,p)=\gcd(m-1,p)=1$. Then $\alpha \in {\Bbb F}_q$. 
\end{lemma}
{\bf Proof}.  From Lemma~\ref{lemma:derivative-of-Sni}, we know that the roots of the derivative of $s_{n,i}(t)$ are the same as those of $s_{n-1,i-1}(t)$. From Lemma~\ref{lemma:Sni-splits-partially}, we know  that under the hypothesis of the present lemma, the roots in question of $s_{n-1,i-1}(t)$ will lie in ${\Bbb F}_{q^{n-1}}$, while those of $s_{n,i}(t)$ will lie in \nq. For a common root then, it must be a root of both of these polynomials, and hence must lie in the intersection of both these fields, which is ${\Bbb F}_q$. \hfill $\Box$

\begin{corollary} \label{corollary:multiple-root-of-Sni}
If  $p$ does not divide both ${n \choose i}$ and ${n-1 \choose i-1}$, then 
there are no non-zero roots of $s_{n,i}(t)$ of multiplicity $m>1$, such that
$\gcd(m,p)=\gcd(m-1,p)=1$.
\end{corollary}
{\bf Proof}. When restricted to ${\Bbb F}_q$, $s_{n,i}(t)={n \choose i}t^i$, and  $s_{n-1,i-1}(t)={n-1 \choose i-1}t^{i-1}$. Thus, if the root is non-zero, the only way both these polynomials can be zero is if the characteristic divides both ${n \choose i}$ and ${n-1 \choose i-1}$. \hfill $\Box$

\begin{lemma} \label{lemma:Sn1-is-permutation}
$s_{n,1}(t)$ is a permutation polynomial over ${\Bbb F}_{q^m}$ if $\gcd(m,n)=1$ and $p$ does not divide $n$.
\end{lemma}
{\bf Proof}. Suppose the pre-image under $s_{n,1}(t)$ of an element comprised of two distinct elements, then, since $s_{n,1}(t)$ is additive, their difference would be its root. But all the roots of $s_{n,1}(t)$ lie in \nq. If $p$ does not divide $n$, then none of these lie in \fq. Thus, under these hypotheses, $s_{n,1}(t)$ would not have any roots in ${\Bbb F}_{q^m}$, where $\gcd(m,n)=1$, and we would get a contradiction.\hfill $\Box$
 
Now we wish to use these notions in the setting of an algebraic function field in one variable $F/K$, where $K = \nq$, and $F$ is an algebraic extension of $K(x)$, where $x$ is transcendental over $K$.

\begin{definition} The field of $(n,q)$-symmetric functions with coefficients in \nq will be denoted
$$  F_s = \nq(s_{n,1}(x), s_{n,2}(x), \ldots, s_{n,n}(x)) \subset K(x),$$
and the field of $(n,q)$-symmetric functions with coefficients in \fq will be denoted
$$  F_s^\phi = \fq(s_{n,1}(x), s_{n,2}(x), \ldots, s_{n,n}(x)) \subset K(x),$$
where the superscript $\phi$ indicates that the values that these functions take on \nq are fixed by $\phi$. Thus, they lie in \fq, as per Lemma~$\ref{lemma:fsphi-in-fq}$. 
\end{definition}

\section{Symmetric extensions of function fields}

Let $F$ and $K$ be as described earlier. Let $E$ be a finite separable extension of $F$, generated by $y$, where $\varphi(y) = 0$, for $\varphi(T)$ an irreducible polynomial in $F[T]$.

In this section we will introduce families of extensions of $F$ whose generators satisfy explicit equations involving only $(n,q)$-symmetric functions. Let $y$ satisfy
                        $$ g(y) = f(x),$$
where $f, g \in F_s^\phi$. If $K=\nq$, this implies that in the residue field of a rational place,  although the class of $x$ and $y$ will assume  values in ${\Bbb F}_{q^n}\cup \infty$, that of $f(x)$ and $g(y)$ will assume values only in ${\Bbb F}_{q}\cup \infty$. Among the Galois extensions that such equations can produce are the two special cases of extensions of Artin-Schreier and Kummer type. In this chapter, we will mainly investigate the case where $f(x)$ is an $(n,q)$-elementary symmetric polynomial. But the techniques to analyse the case of arbitrary $f \in  F_s^\phi$ remain the same and we will consider many such cases in the following chapters.

\subsection{Symmetric extensions of Artin-Schreier type}

In classical Artin-Schreier extensions, the Galois group is a $1$-dimensional vector space over a subfield of the field of constants. We will, however, consider a modified Artin-Schreier extension, where the Galois group is isomorphic to the subgroup of the additive group of  ${\Bbb F}_{q^n}$ comprising elements whose trace in  \fq is zero. 

Before we proceed to describe such modified Artin-Schreier extensions, we reproduce here some definitions and the relations between the different, ramification groups and the genus in extensions of function fields. References are \cite{ParSha1}, \cite{Ser4}, \cite{Sti1} and  \cite{Wei1}.

\begin{definition}
 Let $F'/F$ be a finite separable extension of function fields. Then, the different of $F'/F$, denoted $\diff(F'/F)$, is a divisor in $F'$ given by 
$$\diff(F'/F) = \prod_{P' \in {\Bbb P}(F')}P'^{d(P'|P)},$$ where $d(P'|P)$ is the different exponent of the place $P'$ lying over $P \in \pf$. The degree of $\diff(F'/F)$, denoted $\degdiff(F'/F)$, is  the degree of this divisor. 
\end{definition}

\begin{proposition}[Hurwitz-genus formula]
Let $F'/F$ be a finite separable extension of function fields. Then, 
$$2g(F') - 2 = [F':F](2g(F) - 2) + \degdiff(F'/F).$$
\end{proposition}

\begin{proposition}[Transitivity of different]
Let $F \subseteq F' \subseteq F''$, where $F''/F'$ and $F'/F$ are both finite separable extensions. Let $P \subseteq P' \subseteq P''$, where $P, P'$ and $P''$ are places of $F, F'$ and $F''$, respectively. Then,
$$ \diff(F''/F) = \diff(F''/F') \diff(F'/F). $$
This yields the following relation between the different exponents
$$ d(P''|P) = d(P''|P') + e(P''|P')d(P'|P). $$
\end{proposition}

\begin{definition}[Ramification groups]
Let $F'/F$ be a Galois extension of function fields. Let $P \subseteq P'$, where  $P$ and $P'$ are places in $F$ and $F'$, respectively. The \ith ramification group of $G = \gal(F'/F)$ relative to $P'$ is $G_i = \{s \in G\; | \; s(v) \equiv v \bmod {P'}^{i+1}, \forall v \in {\mathcal O}_{P'} \}$. Then, $G_{-1}$ is the decomposition group, $G_0$ is the inertia group and $G_{-1}/G_0$ is $\gal(K'/K)$, where $K$ and $K'$ are the residue fields of $P$ and $P'$, respectively. Further, $G_i$ is a normal subgroup of $G_{-1}$ as well as $G_{i-1}$, for $i \geq 0$. $G_1$ is a $p$-group and $G_0/G_1$ is a cylic group of order coprime to $p$. 
\end{definition}

\begin{proposition}[Hilbert's different formula]
Let $F'/F$ be a Galois extension of function fields. Let $P \subseteq P'$, where  $P$ and $P'$ are places in $F$ and $F'$, respectively. Then we have that 
$$d(P'|P) = \sum_{i=0}^\infty |G_i| - 1. $$
\end{proposition}

\begin{proposition} \label{lemma:ramification-groups}
Let $F''/F$ be a Galois extension of function fields. Let $H$ be a normal subgroup of $G = \gal(F''/F)$, and $F'$ be its fixed field. Let $P, P'$ and $P''$ be places in $F$, $F'$ and $F''$, respectively, with $P \subset P' \subset P''$. Then the \ith ramification group of $\gal(F''/F')$ relative to $P''$ is $G_i\cap H$. 
\end{proposition}

Now we proceed to describe extensions of modified Artin-Schreier type. 

\begin{proposition} \label{proposition:Art-Sch}
Let $F/K$ be an algebraic function field, where $K=\nq$ is algebraically closed in $F$. Let $w \in F$ and assume that there exists a place $P \in \pf$ such that 
$$ v_P(w) = -m, \, m > 0 \mbox{  and  } \gcd(m,q)  = 1.  $$
Then the polynomial $l(T)-w = a_{n-1}T^{q^{n-1}} + a_{n-2}T^{q^{n-2}}+ \ldots + a_0T - w \in F[T] $ is absolutely irreducible. Further, let $l(T)$ split into linear factors over $K$. Let $E=F(y)$ with
$$  a_{n-1}y^{q^{n-1}} + a_{n-2}y^{q^{n-2}}+ \ldots + a_0y = w. $$
Then the following hold:
\begin{enumerate}
\item $E/F$ is a Galois extension, with degree $[E:F] = q^{n-1}$. $\gal(E/F) = \{\sigma_\beta:y \rightarrow y + \beta\}_{l(\beta)=0}$. 
\item $K$ is algebraically closed in $E$.
\item The place $P$ is totally ramified in $E$. Let the unique place of $E$ that lies above $P$ be $P'$. Then the different exponent $d(P'|P)$ in the extension $E/F$ is given by 
$$ d(P'|P) = (q^{n-1}-1)(m+1).$$
\item Let $R \in \pf$, and $v_R(w) \geq 0$. Then $R$ is unramified in $E$.
\item If $a_{n-1}=\ldots=a_0=1$, and if $Q\in \pf$ is a zero of $w-\gamma$, with $\gamma \in \fq$. Then $Q$ splits completely in $E$.
\end{enumerate}
\end{proposition}
{\bf Proof}. For (i) - (iv), pl. refer \cite{Sul1}. For (v), notice that under the hypotheses, the equation $ T^{q^{n-1}} +T^{q^{n-2}}+ \ldots + T = \gamma $ has $q^{n-1}$ distinct roots in $K$.  \hfill $\Box$

For many of the extensions that we will describe, there exists no place where the hypothesis of Proposition~\ref{proposition:Art-Sch} is satisfied, namely, that the valuation of $w$ at the place is negative and coprime to the characteristic. In particular, we need a criterion for determining the irreducibility of the equations that we will need to use. We provide such a criterion in Proposition~\ref{proposition:irreducibility} and Corollary~\ref{corollary:irreducibility}.

\begin{proposition} \label{proposition:irreducibility}
Let $V$ be a finite subgroup of the additive group of \clo. Then $V$ is a \fp-vector space. Define $L_V(T) = \prod_{v \in V}(T-v)$. Thus, $L_V(T)$ is a separable \fp-linear polynomial whose degree is the cardinality of $V$. Now let $h(T,x) =  L_V(T) - f(x)$, where $f(x) \in \clo[x]$. Then, $h(T,x)=L_V(T) - f(x)$ is reducible over $\clo[T,x]$ iff there exists a polynomial $g(x) \in \clo[x]$ and a proper additive subgroup $W$ of $V$ such that $f(x) = L_{W'}(g(x))$, where $W' = L_W(V)$.
\end{proposition}
{\bf Proof}. Suppose $h(T,x) = g_1(T,x)\ldots g_n(T,x)$ is a factorization of 
$h(T,x)$ over $\clo[T,x]$. From Gauss' lemma, and the fact that $h(T,x)$ is monic in $T$, we see that this is equivalent to a factorization into monic irreducible factors in $\clo(x)[T]$. 

Let $F = \clo(x)$ and view $h(T)$ and $g_i(T)$ as polynomials over $F$. Let $L = F[T]/(g_1(T))$. Then $L$ is a field. Let $y$ be the corresponding root of $g_1$ in $L$.  

Note that $v \in V$ acts on the roots of $h(T)$ over $L$ via $(T \rightarrow T + v)$. Since $y$ is a root of $h(T)$, we see that $h(T)$ splits completely over $L$, into factors $(T-y-v)$, where $v \in V$. 

Also, since $h(T)$ is separable, the factors $g_1(T),\ldots,g_n(T)$ are all distinct. Since the roots of $g_i(T)$ are all conjugate over $F$, it follows that the action of $V$ on the roots of $h(T)$ must respect the $g_i$. Thus, the action of $V$ induces an action on the set $\{g_1(T),\ldots,g_n(T)\}$. 

Let $W$ be the stabilizer of $g_1$ in $V$. Since $V$ is transitive on $\{g_1(T),\ldots,g_n(T)\}$, $|W| = |V/n| = deg(g_1)$. Thus $W$ is isomorphic to the Galois group of $g_1(T)$. It follows that
\begin{eqnarray*}
g_1(T) & = & \prod_{w \in W} (T-y-w), \\
       & = & L_W(T-y), \\
       & = & L_W(T) - L_W(y). 
\end{eqnarray*}
Also note that the factorization of $h(T)$ is nontrivial iff $W$ is a proper
subgroup of $V$.

Now note that the constant term in $g_1(T)$ is $L_W(y)$. Thus $L_W(y) \in F$. In particular, $L_W(y)$ is a polynomial in $x$ with \clo coefficients.  It is useful to note that $L_W : \clo \rightarrow \clo$ is a homomorphism of additive groups with kernel $W$. Finally, let $\hat{W}$ be a complement of $W$ in $V$ and let $W' = L_W(V)$. Then 
\begin{eqnarray*}
h(T) & = & \prod_{v \in \hat{W}}\prod_{w \in W} (T - y - v - w), \\
     & = & \prod_{v \in \hat{W}} L_W(T - y - v), \\
     & = & \prod_{v \in \hat{W}} L_W(T) - L_W(y) - L_W(v), \\
     & = & \prod_{w \in W'}  L_W(T) - L_W(y) - w, \\
     & = & L_{W'}(L_W(T) - L_W(y)), \\
     & = & L_{W'}(L_W(T)) -L_{W'}(L_W(y)), \\
     & = & L_V(T) - L_{W'}(L_W(y)).
\end{eqnarray*}

Comparing this to the equation $h(T) = L_V(T) - f$, we see that $f = L_{W'}(L_W(y))$. This proves the proposition. \hfill $\Box$

The following observation also follows from the proof of Proposition~\ref{proposition:irreducibility}. Let $E = F(y)$ be an extension of $F = \clo(x)$ with Galois group $V$, such that $y$ satisfies the equation $L_V(y) = f(x)$. Let $W$ be a subgroup of index $p$ in $V$, and let $M$ be its fixed field. Then clearly, $z = L_W(y) \in M$. Let $W' = L_W(V) \cong V/W$. Then we have that 
$$ L_{W'}(z) = L_{W'}(L_W(y)) = L_V(y) = f(x). $$  
By noting that the previous equation is of degree $p$, it follows that $M = F(z)$, where $z$ satisfies the equation
$$ L_{W'}(z) = f(x). $$ 

\begin{definition}
For $f(x) \in \clo[x]$, a coprime term of $f$ is a term with non-zero coefficient in $f$ whose degree is coprime to $p$. The coprime degree of $f$ is the degree of the coprime term of $f$ having the largest degree.
\end{definition}

\begin{corollary} \label{corollary:irreducibility}
Let $f(x) \in \clo[x]$. Let there be a coprime term in $f(x)$ of degree $d$, such that there are no terms of degree $dp^i$ for $i>0$ in $f(x)$. Then $L_V(T) -f(x)$ is irreducible for any subgroup $V \subset \clo$.
\end{corollary}
{\bf Proof}. Suppose $f(x)$ is the image of a linear polynomial $\sum a_nx^{p_n}$. Then the coprime term can only occur in the image of the term $a_0x$. But then, the images of the coprime term under $a_nx^{p^n}$, for $n>0$ will have degrees that contradict the hypothesis.

\begin{example}
The equation 
$$ y^{q^2} + y^q + y = x^{q^2 + q} + x^{q^2 + 1} + x^{q+1} $$
is absolutely irreducible over \clo. This follows from Corollary~\ref{corollary:irreducibility} by noting that the coprime degree of the RHS is $q^2+1$, and there are no terms of degree $(q^2 + 1)p^i$, for $i>0$. 
\end{example}

To state the main theorem of this section which describes a general extension of modified Artin-Schreier type using $(n,q)$-elementary symmetric functions, we need  the following lemma.

\begin{lemma} \label{lemma:subextensions}
 Let $F = K(x)$, where $K=\nq, q = p^m, r=m(n-1),$ and $E = F(y)$, where $y$ satisifes the following equation:
$$ y^{q^{n-1}} + y^{q^{n-2}} + \ldots + y = f(x), $$
and $f(x) \in F$ is not the image of any element in $F$ under a linear polynomial. 
Then the following hold:
\begin{enumerate}
\item $E/F$ is a Galois extension of degree $[E:F]=q^{n-1}$. $\gal(E/F) = \{\sigma_\beta: y \rightarrow y + \beta \}_{s_{n,1}(\beta)=0}$ can be identified with the set of elements in \nq whose trace in \fq is zero by $\sigma_\beta \leftrightarrow \beta$. This gives it the structure of a $r$-dimensional \fp vector space.
\item There exists a tower of subextensions 
$$ F=E^0 \subset E^1 \subset \ldots \subset E^{r} = E, $$
such that for $0 \leq i \leq r-1,\; [E^{i+1}: E^i]$ is a Galois extension of degree $p$.
\item Let $\{b_i\}_{1 \leq i \leq r}$ be a \fp-basis for $\gal(E/F)$. Then we can build the tower of subextensions as follows. We set $E^j$ to be  the fixed field of the subgroup of the Galois group that corresponds to the \fp-subspace generated by $\{b_1,b_2,\ldots, b_{r-j}\}$. Then, the generators of $E^{j}$ are $\{y_1,y_2,\ldots,y_j\}$, where $y_1,y_2,\ldots,y_{r}=y$ satisfy the following relations:
\begin{eqnarray*}
       y^p - B_{r}^{p-1} y &=& y_{r-1}, \\
       y_{r-1}^p - B_{r-1}^{p-1} y_{r-1} &=& y_{r-2}, \\
          \vdots &  & \vdots \\
       y_1^p - B_{1}^{p-1}y_1 &=& f(x),
\end{eqnarray*} 
where, 
$$
\begin{array}{rcll}
  \beta_{r,j} &=& b_{r-j+1},  & B_{r} = \beta_{r,r}, \\
  \beta_{r-1,j} &=& \beta_{r,j}^p - B_{r}^{p-1}\beta_{r,j}, & B_{r-1} = \beta_{r-1,r-1},\\
 \vdots & & \vdots &\vdots  \\
  \beta_{1,j}   &=& \beta_{2,j}^p - B_2^{p-1}\beta_{2,j}, &  B_1 = \beta_{1,1}.
\end{array}
$$
\end{enumerate}
\end{lemma}
{\bf Proof}. For (i) refer Proposition~\ref{proposition:Art-Sch} and Proposition~\ref{proposition:irreducibility}. For (ii), note that since  $\gal(E/F)$ is an elementary Abelian group of exponent $p$, we can always find a normal series $\gal(E/K)=G^0 \rhd G^1 \rhd \ldots \rhd G^{r} = 1$ such that $|G^{i+1}/G^i| = p$. Now the existence of a tower $E^0 \subset E^1 \subset \ldots \subset E^{r}$, with $E^{j+1}/E^j$ Galois of degree $p$ is guaranteed by the Fundamental Theorem of Galois Theory by setting $E^i$ to be the fixed field of $G^{i}$. Alternately, we can get a constructive proof for (ii) from the proof for (iii), which follows.
 
For (iii), we set $G^{r-1}$ to be the \fp-span of $b_1$ and $E^{r-1}$ to be its fixed field. Then, the automorphisms of $E^{r}/E^{r-1}$ are given by $y \rightarrow y + ab_1$, with $a \in \fp$. Thus we have that
$$ \prod_{a \in \fp} y - ab_1 = {b_1}^p \prod_{a \in \fp} \frac{y}{b_1} - a  = y^p - {b_1}^{p-1}y = y_{r-1}. $$
Now, we set $G^{r-2}$ to be the \fp-span of $\{b_1, b_2\}$, and iterate this procedure for $E^{r-1}/E^{r-2}$, keeping in mind that the automorphism of $E/F$ given by $y \rightarrow y + b_2$, when pulled down to an automorphism of $E^{r-1}/E^{r-2}$, is given by $y_{r-1} \rightarrow y_{r-1} + {b_2}^p - {b_1}^{p-1}b_2$. By similarly setting $G^{r-j}, j \geq 2$ to be the \fp-span of $\{b_1,b_2,\ldots,b_{j}\}$, and  pulling  down the appropriate automorphisms of $E/F$ to those of  $E^{r-j+1}/E^{r-j}$ we get the other defining equations. The observation made after the proof of Proposition~\ref{proposition:irreducibility} completes the proof. \hfill $\Box$

\begin{corollary} \label{corollary:p-extension}
Let $E,F,K$ be as in Lemma~\ref{lemma:subextensions}. Then, every subextension $E^1$ of $E$ which has degree $p$ over $F$ is of the form $F(z)$, where $z$ satisfies an equation
$$ z^p - Az = f(x),  $$
where $A \in \nq$.
\end{corollary}
It is useful to note that since $E^1$ is the fixed field of a subgroup of $\gal(E/F)$ of index $p$, we can obtain a basis for $\gal(E/F)$ by adding one more element to a basis for this subgroup. 

We are now in a position to state our main theorem for this section.

\begin{theorem} \label{theorem:Art-Sch-symmetric}
Let $F=K(x)$ where $K=\nq, q =p^m \mbox{ and } r=m(n-1)$. Consider the family of extensions $E_i=F(y)$, $2 \leq i \leq n$, where $y$ satisfies the equation
\begin{equation}
y^{q^{n-1}} + y^{q^{n-2}} + \ldots+ y = s_{n,i}(x). \label{equation:Art-Sch-symmetric}
\end{equation}
Then the following hold:
\begin{enumerate}
\item $E_i/F$ is a Galois extension, with degree $[E_i:F] = q^{n-1}$. $\gal(E_i/F)= \{\sigma_\beta:y \rightarrow y + \beta\}_{s_{n,1}(\beta) = 0}.$  
\item The only place of $F$ that is ramified  in $E_i$ is the unique pole  $P_{\infty}$ of $x$. Furthermore, $P_\infty$ is totally ramified in $E_i$. Let $P'_\infty$ denote the unique place of ${\Bbb P}(E_i)$ that lies above $P_\infty$.
\item  Let $m_i$ denote the coprime degree of $s_{n,i}(x)$. We have that
$$ m_i = q^{n-1} + q^{n-2} + \ldots + q^{n-i+1} + 1. $$
The filtration of ramification groups relative to $P_\infty'$ is as follows:
$$ \gal(E_i/F) = G_0 = G_1 = \ldots = G_{m_i+1}, $$
$$ G_{m_i + 2} = \{0\}.$$ 
\item The different exponent $d(P'_\infty|P_\infty)$ is given by 
$$ d(P_\infty|P'_\infty) = (q^{n-1}-1)(m_i+1).$$
\item The genus of $E_i$ is given by 
$$ g(E_i) = \frac{(q^{n-1}-1)(m_i-1)}{2}.$$
\item All other rational places of $F$ split completely in $E_i$ giving
 $$ N(E_i) = q^{2n-1}+1.$$
Thus,  the number of rational places is independent of the choice of $(n,q)$-elementary symmetric polynomial of $x$. 
\end{enumerate}
\end{theorem}                        
{\bf Proof}. (i) follows from Corollary~\ref{corollary:irreducibility} by observing that $p$ times the coprime degree of $s_{n,i}(x)$ is greater than its degree. Let $G = \gal(E_i/F)$. To determine the sequence of ramification groups of $G$ relative to $P_\infty'$, we will use the function $i_G$ defined on the elements of the $G$ as follows: 
$$ i_G(s) \geq k+1 \Leftrightarrow s \in G_k.$$

For every subgroup $H$ of $G$, define $(G/H)_k$ to be the $k^{th}$ ramification group of the fixed field of $H$, relative to its unique place lying above $P_\infty$. Then, define the function $i_{G/H}$ on the elements of $G/H$ as follows:
$$ i_{G/H}(\overline{s}) \geq k+1 \Leftrightarrow \overline{s} \in (G/H)_k. $$

Note that for all subgroups $H$ of index $p$ in $G$, we know that
\begin{eqnarray*}
i_{G/H}(\overline{s}) &=& \left\{ \begin{array}{ll} \infty & \mbox{if } \overline{s} = 0, \\
					m_i + 2 & \mbox{else}. 
\end{array}
\right.
\end{eqnarray*}

We claim that if $K$ is any proper subgroup of $G$ then the average of the values of the function $i_{G/K}$ on the non-zero elements of $G/K$,
$$ \displaystyle{\mbox{avg}}_{s\in G/K, s \neq 0} i_{G/K}(s) = m_i + 2. $$

To see this, suppose that $|G/K| = p^l$. Consider the subgroups $H$ such that $K \subset H \subset G$, with $[G:H]=p$. These are in $1:1$ correspondence with the subgroups of index $p$ in $G/K$, which number $\frac{p^l - 1}{p-1}$. A non-zero element $s \in G/K$ is contained in exactly $\frac{p^{l-1}-1}{p-1}$ of these. Now, from \cite{Ser4}, Ch. IV, Proposition 3, we get
$$ i_{G/H}(\overline{s}) = \frac{1}{p^{l-1}} \sum_{s \rightarrow \overline{s}} i_{G/K}(s). $$
Since each non-zero $s \in G/K$ is nontrivial in exactly $p^{l-1}$ of the $G/H$, so that summing over all non-zero $s \in G/K$, we get
$$ \sum_{K\subset H \subset G, [G:H]=p \atop \overline{s} \in G/H, \overline{s} \neq 0} i_{G/H}(\overline{s}) = \sum _{s\in G/K} i_{G/K}(s). $$ 

But as noted previously, each $i_{G/H}(\overline{s}) = m_i +2$. Since each side in the previous equation has the same number of terms, the average of the RHS = average of the LHS = $m_i + 2$. This proves the claim.

It now follows that $i_G(s) = m_i + 2$ for all non-zero $s \in G$. For
suppose not. Then there must exist $s \in G$ such that $i_G(s) > m_i + 2$. Since $i_G$ is constant on cyclic subgroups, we see that $\langle{s}\rangle \neq G$. But then the average on $G/\langle{s}\rangle$ will be less than $m_i + 2$, giving a contradiction.

(iv) now follows from Hilbert's different formula. (v) is a straightforward application of the Hurwitz-genus formula. (vi) follows from Proposition~\ref{proposition:Art-Sch}. We obtain the number of rational places in $E_i$ as follows. Each of the $q^n$ finite places in $F$ splits completely in $E_i$, giving $(q^n)(q^{n-1}) = q^{2n-1}$ rational places in $E_i$. Also, $P_\infty$ ramifies totally in $E_i$, giving a sum total of $q^{2n-1}+1$ rational places in $E_i$. \hfill $\Box$

\begin{example} \label{example:n=3,i=2}\rm (n=3, i=2) Let $F=K(x)$, where $K = {\Bbb F}_{q^3}$. Let $E=F(y)$ where $y$ satisfies the equation
$$ y^{q^2} + y^q + y = x^{q^2 + q} + x^{q^2 + 1} + x^{q+1}.$$
All rational places of $F$, except $P_\infty$,  split completely in $E$, giving a total of $(q^3)(q^2) = q^5$ rational places. Also, $P_\infty$ ramifies totally in $E$ giving one rational place. Thus, $N(E) = q^5+1$ rational places. The genus $g(E) =  \frac{(q^2-1)(q^2)}{2}$. In a later section we will see that this extension attains the Oesterle lower bound on genus for $q=2$, \ie, over ${\Bbb F}_{q^3} = {\Bbb F}_8$.
\end{example}

\begin{example} \label{example:n=3,i=2,subfield}\rm (n=3, i=2) In this example, we discuss a subfield of the Example~\ref{example:n=3,i=2}. Let $E,F,K$ all be as in Example~\ref{example:n=3,i=2}. Then consider the extension $E^1/F$ where $E^1=F(y_1)$, and $y_1$ satisfies the equation
$$ y_1^q + (1+b^{q^2-q})y_1 = x^{q^2 + q} +  x^{q^2 + 1} + x^{q+1},$$ 
where $b(\neq 0) \in K$ is an element whose trace in \fq is zero. Then $E^1$ is a subfield of $E$ and $E = E^1(y)$, with $y$ satisfying
$$ y^q - b^{q-1}y = y_1. $$  
$P_\infty$ is totally ramified in $E^1$ and all other rational places of $F$ split completely in $E^1$ (This is clear since $E^1$ is a subfield of $E$, in which these statements are true). Let $P^1_{\infty}$ be the unique place of $E^1$ lying above $P_\infty$. Then we have that $d(P^1_{\infty}|P) = (q-1)(q^2 + 2) $ and $g(E^1) = \frac{(q-1)(q^2)}{2}$. Also, $N(E^1) = q^4+1$. $E^1$ also attains the Oesterle lower bound on genus for $q=2$.  
\end{example}

\begin{example} \rm (n=4, i=3) Let $F=K(x)$, where $K = {\Bbb F}_{q^4}$. Let  $E=F(y)$ where $y$ satisfies the equation
$$ y^{q^3} + y^{q^2} + y^q + y = x^{q^3 + q^2 + q} + x^{q^3 + q^2 + 1} + 
x^{q^3+q+1} + x^{q^2 + q + 1}. $$
All rational places of $F$, except $P_\infty$,  split completely in $E$. $P_\infty$ ramifies totally in $E$. This gives us $N(E)=q^7+1$. The genus $g(E) =  \frac{(q^3-1)(q^3+q^2)}{2}$.
\end{example}

\begin{definition}
A function field $F/K$, where $K = \nq$ is said to be median if $N(F) = q^n + 1$. Thus, the number of its rational places is exactly in the middle of the range allowed by the Hasse-Weil bound. 
\end{definition}

\begin{proposition}
Let $F=K(x)$, where $K = {\Bbb F}_{q^m}, \; \gcd(m,n) = \gcd(p,n) =1$. Let $E_i = F(y)$ where $y$ satisfies the equation 
$$ y^{q^{n-1}} + y^{q^{n-2}} + \ldots + y = s_{n,i}(x). $$
Then $E_i$ is median\footnote{Function fields which are median over infinitely many extensions of their field of constants (\fq in our case) are called {\em exceptional}. There is a rich theory to such function fields. For instance, it is known \cite{Fri1} that the roots of their zeta-function occur in cliques as roots of unity times each other.}. It retains this property if we replace $s_{n,i}(x)$ with any other polynomial such that the resulting equation is irreducible.
\end{proposition}
{\bf Proof}. From Lemma~\ref{lemma:Sn1-is-permutation}, we know that under the hypothesis,  $ y^{q^{n-1}} + y^{q^{n-2}} + \ldots+ y$ simply permutes the elements of ${\Bbb F}_{q^m}$.  Then, for every value that $x$ can take in ${\Bbb F}_{q^m}$, we have exactly one solution for $y$.  \hfill $\Box$

\subsubsection{Special case: $i=n$ (trace-norm)}

It has been known that in extensions of the form $E/F$ where $E=F(y)$ with
$$  s_{n,1}(y) = s_{n,n}(x),  $$ 
all the rational places of $F$, except $P_\infty$, split completely, yielding many rational places in $E$. Extensions of this form have been referred to as ``trace-norm'' extensions. However, we will treat this type of extension as a special case of the generalized symmetric extensions of Theorem~\ref{theorem:Art-Sch-symmetric}. In the notation of Theorem~\ref{theorem:Art-Sch-symmetric}, the trace-norm extension is $E_n$. 

The most famous example of a trace-norm construction is the Hermitian function field with full field of constants ${\Bbb F}_{q^2}$ (\ie, the case $n = 2$), where the trace and norm are taken down to ${\Bbb F}_q$. This is a maximal function field in the Hasse-Weil sense. Whenever $q \neq p$, we may take the trace and norm down to a smaller subfield of ${\Bbb F}_{q}$. All three - the degree of the extension, the number of rational places, and the genus - for such a construction increase as we take traces and norms to smaller  subfields. The maximum for each of these is attained when we take the trace and norm down to the prime field ${\Bbb F}_p$.

For a function field over $K$, we may then construct trace-norm extensions by taking these into any subfield of $K$. In the language of $(n/m,q^m)$-elementary symmetric polynomials, the most general form of the trace-norm extension is given by 
$$ s_{\frac{n}{m},1}(y) = s_{\frac{n}{m},\frac{n}{m}}(x), $$
where the trace and norm are being taken in the subfield ${\Bbb F}_{q^m}$.
Now we would like to see how the ratio $N/g$ in such extensions varies as we vary $m$ from $1$ to the value of the greatest proper divisor of $n$. In other words, we would like to see how this ratio varies as we change the subfield of $K$ in which we take the trace and norm.

\begin{lemma} \label{lemma:N/g-trace-norm}
Let $F=K(x)$ where $K=\nq$. Let $m \neq n, m|n$ and $\frac{n}{m} = r$. Let $E=F(y)$, where $y$ satisfies the equation
\begin{equation}
y^{q^{m(r-1)}} + y^{q^{m(r-2)}} + \ldots + y = x^{q^{m(r-1)}+ \ldots + q^m + 1}. \label{equation:N/g-trace-norm}
\end{equation}
Thus, the trace and norm are being taken to the field ${\Bbb F}_{q^m}$. Then the following holds:
\begin{enumerate}
\item The ratio $N/g$, of the number of rational places to genus, decreases with decreasing $m$.
\item If $n \equiv 0\bmod 2$, the maximum value of  $N/g$ for the extension of this type is obtained by taking trace and norm down to the field of cardinality $q^{\frac{n}{2}}$.
\end{enumerate}
\end{lemma}
{\bf Proof}. In the general case, $N(E) = q^{2n-m}+1$, and $g(E) = \frac{q^m(q^{n-m} - 1)^2}{2(q^m-1)}$. Thus, 
\begin{eqnarray*}
 \frac{N(E)}{2g(E)} & = &  \frac{(q^m-1)(q^m+q^{2n})}{(q^n - q^m)^2}.
\end{eqnarray*}
The numerator increases with increasing $m$ while the denominator decreases with increasing $m$. The result follows.   \hfill $\Box$

\noindent{{\bf Note}: Lemma~\ref{lemma:N/g-trace-norm}, (ii) is the case of the Hermitian function field.}

\begin{corollary}
For extensions of the form given by Lemma$~\ref{lemma:N/g-trace-norm}$, the lowest value of the ratio $N/g$ is reached when we take norms and traces down to ${\Bbb F}_p$.
\end{corollary}

There are other extensions of function fields using trace and norm, and one of them is given below. 

\begin{example} \rm Let $F=K(x)$, where $K={\Bbb F}_{q^2}$. Let $E=F(y)$ with
$$ y^q + y = \frac{x^{q+1}}{x^q + x}. $$
This is the function field at the second step of the tower of function fields attaining the Drinfeld-Vladut bound from \cite{GarSti2}. The only places of \pf that are ramified in $E$ are $P_{\infty}$ and $\{P_{\alpha}\}_{s_{2,1}(\alpha) = 0, \alpha \neq 0}$. These are totally ramified, each with different exponent $2(q-1)$. All other rational places split completely. Thus we have that $N(E) = q^3-q^2+2q$ and $g(E) = (q-1)^2$.
\end{example}

\begin{lemma} 
The extension $E_n$ described in Theorem~$\ref{theorem:Art-Sch-symmetric}$ attains the Hasse-Weil bound for $n=2$, for all values of $q$. For $n>2$, it does not attain the Hasse-Weil bound for any value of $q$. 
\end{lemma}
{\bf Proof}. Observe that for $n>2$, 
$$  \frac{q{(q^{n-1} - 1)}^2}{2(q-1)} > \frac{q^{\frac{n}{2}}(q^{\frac{n}{2}} -
 1)}{2}. $$
The lemma then follows from a well-known result that says that a function field over ${\Bbb F}_l$ cannot attain the Hasse-Weil bound for genus $g>\frac{\sqrt{l}(\sqrt{l}-1)}{2}$, cf. \cite{Sti1}, Ch. V.3.3. \hfill $\Box$

\subsection{Symmetric extensions of Kummer type}                                                                       
We now study extensions whose Galois group is a subgroup of the 
multiplicative group $K^*$. For this we will need that the field 
contain a primitive \jth root of unity $\xi_j$ for some $j$ coprime to $p$. In particular we know that $K$ contains $\xi_j$ for $j = \frac{q^n-1}{q-1}$.                                 

\begin{theorem} \label{theorem:Kummer-symmetric}
Let $F=K(x)$ where $K=\nq$. Let  $E_{i,Kum}=F(y)$, $1\leq i \leq n-1$, where $y$ satisfies the equation
\begin{equation}
 y^{\frac{q^n-1}{q-1}} = s_{n,i}(x). \label{equation:Kummer-symmetric}
\end{equation}
Then the following hold:
\begin{enumerate}
\item $E_{i,Kum}/F$ is a cyclic Galois extension, with degree $[E_{i,Kum}:F] = \frac{q^n-1}{q-1}$. $\gal(E_{i,Kum}/F) = \{ \sigma_j: y \rightarrow y{\xi}^k \}_{1 \leq k \leq \frac{q^n-1}{q-1}}.$
\item The only places of $F$ that are ramified in $E_{i,Kum}$ are $P_\infty$ and $\{P_\alpha\}_{s_{n,i}(\alpha)=0}$. Let $v_P = v_P(s_{n,i}(x))$.  Define $r_\alpha = \gcd([E_{i,Kum}:F],v_{P_\alpha}) > 0$ and $r_{\infty}= \gcd([E_{i,Kum}:F],v_{P_\infty}) > 0$. Then we have that
\begin{eqnarray*}
 e(P'_\alpha|P_\alpha) &=& \frac{[E_{i,Kum}:F]}{r_\alpha}, \\
 e(P'_\infty|P_\infty) &=& \frac{[E_{i,Kum}:F]}{r_{\infty}}.
\end{eqnarray*} 
Since the extension is tame, the different exponents are given by
 $$d(P'|P) = e(P'|P) - 1,\; \forall P' \in \pe.$$
Also, $v_{P_0}(s_{n,i}(x))= \frac{q^i-1}{q-1}$ and $v_{P_\infty}(s_{n,i}(x))= q^{n-i}(\frac{q^i-1}{q-1})$. 
\item All other rational places of ${\Bbb F}_{q^n}(x)$ split completely in $E_{i,Kum}.$
\end{enumerate}
\end{theorem}
{\bf Proof}. The proofs for (i) and (ii) will need standard results on Kummer extensions, cf. \cite{Sti1}, Ch. III.7.3. Also note that $r_\alpha, r_\infty < \frac{q^n-1}{q-1}$. For (iii) we note that from Lemma~\ref{lemma:fsphi-in-fq}, $s_{n,i}(\gamma) \in \fq, \forall \gamma \in \nq$, and therefore, it has $\frac{q^n-1}{q-1}$ pre-images under the norm map. \hfill $\Box$

\begin{lemma}
$$E_{i,Kum} \cong E_{n-i,Kum}.$$
\end{lemma}
{\bf Proof}. By making the substitution $y = xy'$, and then using Lemma~\ref{lemma:Sni-Sn(n-i)}. \hfill $\Box$

This immediately leads to the following corollary.
\begin{corollary} \label{corollary:isomorphic-to-trace-norm}
The extensions $E_{1,Kum}$ and $E_{n-1,Kum}$ of Theorem~$\ref{theorem:Kummer-symmetric}$ are both isomorphic to the trace-norm extension $E_n$ of Theorem~$\ref{theorem:Art-Sch-symmetric}$. 
\end{corollary}

\begin{example} \rm (n=3, i=2)  Let $F=K(x)$, where $K = {\Bbb F}_{8}$. Let $E=F(y)$ where $y$ satisifes the equation
$$ y^7 = x^6 + x^5 + x^3 \; ( = s_{3,2}(x)). $$
$s_{3,2}(x)$ has three distinct zeros, other than $0$ itself, which has multiplicity three. $N(E) = 33$ and $g(E) = 9$. From Corollary~\ref{corollary:isomorphic-to-trace-norm}, we can see why these are the same values as the corresponding trace-norm extension. 
\end{example}

There are other examples of Kummer extensions using $(n,q)$-elementary symmetric polynomials alone. An example from \cite{GarSti3} is given below. 

\begin{example} \rm
Let $F = \fq(x)$, $q=p^e, e > 1, m = \frac{q-1}{p-1}$, and let $E=F(y)$ with
$$ y^m = (1+x)^m + 1.$$
Here $m$ rational places in \pf ramify, while all others split completely in $E$. Notice that $P_\infty$ splits completely in this extension, unlike the extensions of Theorem~\ref{theorem:Kummer-symmetric}. This construction, when iterated, gives an asymptotically good tower \cite{GarSti3}.
\end{example}

\section{Number of rational places versus genus}

In order to arrive at the Oesterle lower bound on genus $g$ of a function field for a specified field of constants \fq and specified number of rational places $N = L +1$, one must go through an algorithm of sorts, that is given below \cite{Sch2}.
\renewcommand{\labelenumi}{\rm (\arabic{enumi})}
\begin{enumerate}
\item Let $m$ be the unique integer for which 
$$ {\sqrt{q}}^{\,m} < L \leq {\sqrt{q}}^{\,m+1}. $$
\item  Define 
$$ u=\frac{{\sqrt{q}}^{\,m+1} - L}{L{\sqrt{q}}-{\sqrt{q}}^{\,m}} \in [0,1).$$ 
\item  Let $\theta_0$ be the unique solution of the trigonometric equation
$$ cos{\frac{m+1}{2}\theta} + u\, cos{\frac{m-1}{2}\theta} = 0 $$
in the interval $[\frac{\pi}{m+1}, \frac{\pi}{m}]$. \\
\item  Then we have that 
$$g \geq \frac{(L-1)\sqrt{q}\, cos \theta_0 + q - L}{q+1-2\,\sqrt{q}\, cos \theta_0}.$$
\end{enumerate}
\renewcommand{\labelenumi}{\rm (\roman{enumi})}

\begin{example} \rm
If $q=8$ and $N=17$, then we have that $m=2$, $u=0.1779$ and $\theta_0 =1.1472$ yielding $g \geq 1.414$. Thus for a function field  over  ${\Bbb F}_{8}$ to have $17$ rational places, it must have genus at least 2.
\end{example}

\begin{example} \rm
If $q=8$, and $N=33$, then we have that $m=3$, $u=0.47$ and $\theta_0 = 0.88735$, yielding $g \geq 5.779$. Thus for a function field over  ${\Bbb F}_{8}$ to have $33$ rational places, it must have genus at least 6.
\end{example}

\begin{example} \rm
If $q=16$, and $N=129$, then we have that $m=3, u=0.2857$ and $ 
\theta_0=0.87752$, yielding $g \geq 17.88$. Thus for a function field over ${\Bbb F}_{16}$ to have 129 rational places, it must have genus at least 18.
\end{example}

\begin{example} \rm
If $q=27$, and $N=244$, then we have that $m=3, u=0.433$ and $ 
\theta_0=0.90754$, yielding $g \geq 25.16$. Thus for a function field over ${\Bbb F}_{27}$ to have 244 rational places, it must have genus at least 26.
\end{example}

Now we investigate the performance of symmetric function fields with respect to the known (Hasse-Weil, Oesterle) bounds on the number of rational places.

Theorem~\ref{theorem:Art-Sch-symmetric} describes, for a specific value of $n$, $n-1$ different symmetric extensions of the rational function field corresponding to  $i=2,3,\ldots,n$. The case of $i=n$ is the trace-norm extension. We now compare these $n-1$ different extensions for various values of $q$.
\begin{center}
\begin{tabular}{|c|c|c|c|c|c|c|c|} \hline
$n$ & $q$ & $N(E)$  & $g_n(E)$ & $g_{n-1}(E)$ & $g_{n-2}(E)$ & $g_{n-3}(E)$ & Oesterle  \\ \hline
3 & 2 & (17)  &  	& $(2)^*$ &	&	  			&$2^*$     	\\ 
3 & 2 & 33    & 9 	& $6^*$ &	&	  			&$6^*$     	\\ 
3 & 4 & 1025  & 150 	& 120   &	&		    	    	& 74       	\\
3 & 8 & 32769 & 2268 	& 2016  & 	&				& 903      	\\
3 & 3 & 244   & 48   	&  36   & 	&			        & 26       	\\
3 & 9 & 59050 & 3600	& 3240  &	&				& 1374     	\\
3 & 5 & 3126  & 360	& 300   &	&				& 167      	\\
3 & 7 & 16808 & 1344	& 1176  &	&				& 560      	\\
3 & 11& 161052& 7920	& 7620  &	&				& 2808     	\\
4 & 2 & 129   & 49	& 42	& 28	&				& 18 	   	\\
4 & 4 & 16385 & 2646	& 2560  & 2016	&				& 667	   	\\
4 & 3 & 2188  & 507	& 468   & 351	&				& 152	   	\\
4 & 5 & 78126 & 9610	& 9300	& 7750	&				& 2071		\\
5 & 2 & 513   & 225	& 210 	& 180	& 120 				& 57 		\\
\hline
\end{tabular}
\end{center}
{\bf Table 1}: {\footnotesize Number of rational points and genus for some members of the  families of symmetric function fields. The entry $g_j(E)$ denotes the genus of the extension described by Theorem~\ref{theorem:Art-Sch-symmetric} obtained for $i=j$. Note that $g_n(E)$ is the genus of the trace-norm extension. * indicates that the function field attains the Oesterle bound, and parentheses indicate a subfield of the symmetric function field.}

As  Theorem~\ref{theorem:Art-Sch-symmetric} tells us, the extension $E_2$ always has the lowest genus, while the extension $E_n$ (the trace-norm extension) has the highest genus. In other words, for purposes of the ratio of $N/g$, the trace-norm extension is actually the {\em worst} of the symmetric extensions, while the extension that uses the second $(n,q)$-symmetric polynomial is the {\em best}. 

We now look more closely at the case of $n=3,\;i=2$, since this contains a maximal member in the the Oesterle sense. Specifically, we provide a uniformizing parameter at the place $P_\infty'$.
\begin{example} \rm \label{example:n=3,i=2,contd}
Let $F=K(x)$, where $K = {\Bbb F}_{q^3}$. Let $E=F(y)$ where $y$ satisfies the 
equation
$$ y^{q^2} + y^q + y \;=\; x^{q^2+q} + x^{q^2+1} + x^{q+1}. $$
Except for $P_{\infty}$, which is totally ramified, all other rational places in $F$ split completely in $E$, giving a total of $q^5+1$ rational places. A uniformizing parameter for $P_\infty'$ can be obtained as follows. First we make the substitutions $x = X/Z \mbox{ and } y = Y/Z$ to get the following homogeneous equation:

$$Y^{q^2}Z^q + Y^qZ^{q^2} + YZ^{q^2+q-1} \;=\; X^{q^2+q}+X^{q^2+1}Z^{q-1}+X^{q+1}Z^{q^2-1}, $$
where, we look at $Y=1, Z = 0$ (\ie, $P'_\infty$), and we get
$$Z^q + Z^{q^2} + Z^{q^2+q-1} = X^{q^2 + q} + X^{q^2+1}Z^{q-1} + X^{q+1}Z^{q^2-1}. $$
Note that the valuations for $x,y,X,Y \mbox{ and } Z$ at $P'_\infty$ are $-q^2, -(q^2+q), q, 0 \mbox{ and } q^2+q$, respectively. Now, we expand $Z$ in a power series in $X$:
$$ Z = X^{q+1} + X^{2q} + \frac{X^{2q+1}}{\pi}, $$
where $\pi$ has valuation of the $q^{th}$ root of $X$, \ie, it is a uniformizing parameter for $P'_\infty$. Thus we get a uniformizing parameter $\pi$ for $P'_\infty$ given by: 
$$ \pi = \frac{X^{2q+1}}{Z - X^{q+1} + X^{2q}} = \frac{x^{2q+1}}{y^{2q} - x^{q+1}y^q - x^{2q}y}. $$
Now, let $\sigma_\beta \in \gal(E/F)$ be such that $\sigma_\beta: y \rightarrow y + \beta$, where $\beta$ is a non-zero element of ${\Bbb F}_{q^3}$ whose trace in \fq is zero. Then,
$$ \sigma_\beta(\pi) - \pi = \frac{x^{2q+1}[\beta^{2q} + 2\beta^q y^q - \beta^q x^{q+1} - \beta x^{2q}]}{[(y+\beta)^{2q} - x^{q+1}(y+\beta)^q - x^{2q}(y+\beta)][y^{2q} - x^{q+1}y^q - x^{2q}y]}. $$
By noting that the valuation of both the terms in the denominator is the same, since they are conjugates, and then observing that the term in the numerator whose valuation will dominate is the last one, we get the valuation of $\sigma_\beta(\pi) - \pi$, which is independent of the choice of $\beta$:
$$ v_{P'_\infty}(\sigma_\beta(\pi) - \pi) = q^2 + 2. $$
This again gives us the sequence of ramification groups:
$$ G_0 = G_1 = \ldots = G_{q^2 +1}, $$
$$ G_{q^2 + 2} = \{0\}. $$
Thus the different exponent and the degree of the different for the extension $E/F$ are 
$$ d(P_\infty'|P_\infty) =  \degdiff(E/F) = (q^2 + 2)(q^2 -1).  $$
This agrees with the results given in Theorem~\ref{theorem:Art-Sch-symmetric}.

\end{example}

We now wish to compare extensions of the rational function field having the typical (to our family of extensions) jumps in the sequence of ramification groups relative to the unique ramified place but having different degrees as extensions of the rational function field. 
                                                                               \begin{theorem}
Let $F=K(x)$ and $K = \fq$. Let $E/F$ be an extension of function fields of degree $d$, where the sequence of ramification groups relative to the unique ramified place is of the form 
$$ G_0=G_1=\ldots=G_k, $$
$$ G_{k+1} = \{0\}. $$ 
Further let all other places in $F$ split completely in $E$. Then the ratio of $N/g$ for $E$ decreases with increasing $d$.
\end{theorem}
{\bf Proof}. We have that 
\begin{eqnarray*}
N &=& dq + 1, \\2g & = & 2 -2d + (d-1)(k+1) = (d-1)(k-1), \\
\frac{N}{2g} &=& \frac{dq + 1}{(d-1)(k-1)}.
\end{eqnarray*}
Then differentiating w.r.t $d$, we get that 
$$ \left[\frac{N}{2g}\right]' = \frac{(q+1)(1-k)}{[(d-1)(k-1)]^2}. $$
Now notice that $k \geq 2$ for $g(E) > 0$. But in that case the numerator is always nonpositive. Hence the result. \hfill $\Box$

Thus, we get higher $N/g$ ratios by looking at subfields that have identical jumps in their sequence of ramification groups. The smaller the subfield, the higher is the $N/g$ ratio. In a sense, the rational function field is the limiting case of this behaviour, with $N/g = \infty$. 

Thus ray class field extensions, which are maximal Abelian extensions with a certain property, have lower $N/g$ ratios than other Abelian extensions having identical jumps in their sequence of ramification groups.

\begin{example} \rm
 Consider the extension of Example~$\ref{example:n=3,i=2,contd}$. Let $q=2$ and  $G= \gal(E/F)$. Then the filtration of ramification groups at $P'_\infty$ is $G = G_0 = G_1 = \ldots = G_5$ and $G_6$ is trivial. In the ray class field constructions of function fields with many rational places \cite{Lau1}, there is a function field over ${\Bbb F}_8$ with an identical sequence of ramification groups. It has degree $8$, as an extension of the rational function field,  and genus $g=14$, for $N=65$ rational places. Thus,  the ratio  $N/g = 4.64$. The function field in our example has $N=33$ and $g=6$, giving a ratio of $N/g= 5.5$.
\end{example}

\section{Generalization of the Hermitian function field}

 The simplest example of a symmetric function field is the Hermitian function field itself, described below. 

Let $F={\Bbb F}_{q^2}(x)$ and $E=F(y)$ where $y$ satisfies the equation 
$$ y^q + y \,=\, x^{q+1}. $$
Then $E$ is called the Hermitian function field. We have that $ g(E) = \frac{q(q-1)}{2}$, and $ N(E) = q^3 + 1$, which is  the maximum number allowed by the Hasse-Weil bound for this value of genus.  Furthermore, the Hermitian function field is the unique maximal function field over ${\Bbb F}_{q^2}$ of genus $g \geq \frac{q(q-1)}{2}$ \cite{FuhSti1}.  The Hermitian function field remains maximal for all values of $q$. Other symmetric function fields do not exhibit such a uniform performance. For instance, the function fields of  Example~\ref{example:n=3,i=2} and Example~\ref{example:n=3,i=2,subfield} attain the the Oesterle bound for $q=2$, while for higher values of $q$, they deviate considerably from the Oesterle bound. The other unique feature of the Hermitian function field is the extremely large size of its automorphism group \cite{Sti3}.

It has often been implicitly assumed that function field $E_n$ described by Theorem~\ref{theorem:Art-Sch-symmetric}, called the trace-norm function field, is the generalization of the Hermitian function field for $n \geq 3$. However, we argue below that it is the function field $E_2$ described by Theorem~\ref{theorem:Art-Sch-symmetric} that is more the generalization. It is clear that both coincide for $n=2$ and in this case, they are just the Hermitian function field. For $n \geq 3$, we must decide which is the more appropriate generalization of the Hermitian function field.  We now show the following similarities.

\vspace{0.2in}
\noindent{\bf A. Genus}\\
We have already observed that for  the $n-1$ families of symmetric extensions $E_i, 2  \leq i \leq n$ described by Theorem~\ref{theorem:Art-Sch-symmetric},
$$ g(E_2) < g(E_3) \ldots < g(E_n). $$
Thus, among these families $E_2$ has the lowest genus while the trace-norm function field has the highest.

\vspace{0.2in}
\noindent{\bf B. Automorphisms}\\
It is well known \cite{Sti2} that for every $\delta$ and $\tau$ in ${\Bbb F}_{q^2}$ that satisfy $\tau^q + \tau = \delta^{q+1}$, there is an automorphism $\sigma$ of the Hermitian function field  given by
\begin{eqnarray*}
\sigma(x)&=& x+ \tau, \\
\sigma(y)&=& y+ x \tau^q + \delta.
\end{eqnarray*}

We can state the general theorem about automorphisms of the symmetric function field $E=F(y)$ where $y$ satisfies (\ref{equation:Art-Sch-symmetric})  for $i=2$.
\begin{theorem} \label{theorem:generalized-Hermitian}
Let $F=K(x)$ where $K=\nq$. Let $E=F(y)$, where $y$ satisfies the equation
\begin{equation}
 y^{q^{n-1}} + y^{q^{n-2}} + \ldots+ y = s_{n,2}(x). \label{equation:generalized-Hermitian} 
\end{equation}
Let  $(\delta, \tau)$ satisfy $\delta, \tau \in K$ and 
$$  \delta^{q^{n-1}} + \delta^{q^{n-2}} + \ldots+ \delta = s_{n,2}(\tau). $$ 
Let $ m = \lfloor \frac{n-1}{2} \rfloor $. Then there exists an automorphism $\sigma$ of $E$ given by
$$ \begin{array}{lcl}
 \sigma(x) & = & x+\tau,  \\
 \sigma(y) & = & y + x\tau^q + x\tau^{q^2} + \ldots + x\tau^{q^{n-1}} + \delta.
\end{array} $$
The set of these automorphisms keeps $F$ fixed and forms a subgroup
$\Gamma$ of order $q^{2n-1}$ of the full automorphism group of the function field $E/K$. $\Gamma$ acts transitively on the set of finite rational places of $E$.
\end{theorem}
{\bf Proof}. On expanding $s_{n,2}(x + \tau)$, we observe that there are, apart from all the terms of $s_{n,2}(x)$ and $s_{n,2}(\tau)$, cross terms in $x$ and $\tau$ with all possible pairs of exponents $\{q^i,q^j\}$ for $x$ and $\tau$, where $ i \neq j,\;0 \leq i,j \leq n-1 $. Thus, there are $n(n-1)$ such cross terms, all distinct. On expanding 
$[\sigma(y)]^{q^{n-1}} + [\sigma(y)]^{q^{n-2}} + \ldots + \sigma(y)$, we get
$n(n-1)$ cross terms in $x$ and $\tau$, all distinct, with the same exponent pairs as earlier. \hfill $\Box$

\vspace{0.2in}
\noindent{\bf C. Places of degree two}\\
Let us denote by $N_m(E)$ the number of rational places of degree $m$ in a function field $E$. Thus, $N_m$ is the number of Galois conjugacy classes of ``new" rational points of the curve over the extension of degree $m$ of the original field of definition (\ie, points that belong strictly to the extension, and not to the original field).  For the Hermitian function field over ${\Bbb F}_{q^2}$, $N_2 = 0$. This means that there are no ``new" rational places over  ${\Bbb F}_{q^4}$ --  the number of rational points over both the fields is $q^3 +1$. The same phenomenon occurs in each of the other ``Deligne-Lusztig'' curves, namely the Suzuki and the Ree curves.
 
\begin{lemma} \label{lemma:N_2=0}
Let the hypotheses be as in Theorem~$\ref{theorem:Art-Sch-symmetric}$.
Then  $N_2(E_i) \neq 0$ only if there exist ${\alpha, \beta}$, such that
$$ \beta^{q^{n-1}} + \beta^{q^{n-2}} + \ldots+ \beta = s_{n,i}(\alpha), $$
with  $\alpha \in {\Bbb F}_{q^{2n}} \setminus {\Bbb F}_{q^n}$ such that  
$s_{n-1,i-1}(\alpha) \in {\Bbb F}_{q^n}$.
\end{lemma}
{\bf Proof}. Let $(x,y) = (\alpha,\beta)$ satisfy (\ref{equation:Art-Sch-symmetric}). Then, raising the equation to the $q^{th}$ power  and subtracting the original from it, we get 
$$ \beta^{q^n} - \beta = [s_{n-1,i-1}(\alpha)]^q(\alpha^{q^n} - \alpha). $$
Now raising both sides to the power of $q^n$, we get that if $\alpha$ and $\beta$ are 
elements of ${\Bbb F}_{q^{2n}}$, then
$$ \beta - \beta^{q^n} = [s_{n-1,i-1}(\alpha)]^{q^{n+1}}(\alpha - \alpha^{q^n}).$$
Then adding the two equations, we get that either $\alpha=\alpha^{q^n}$ or 
$[s_{n-1,i-1}(\alpha)]^q=[s_{n-1,i-1}(\alpha)]^{q^{n+1}}]$. But since $\alpha \notin \nq$, and \nq is closed under taking $q^{th}$ powers, this gives us the result. \hfill $\Box$

\begin{example}
In the case of the Hermitian function field, since $s_{1,1}(x)=x$, Lemma~$\ref{lemma:N_2=0}$ precludes any possibility for $N_2 \neq 0$. 
\end{example}

We now need another lemma.

\begin{lemma} \label{lemma:neq-2-3-6}
For $n \neq 3,4,6$, There exists no $\alpha \in {\Bbb F}_{q^{2n}} \setminus {\Bbb F}_{q^n}$ such that $s_{n-1,1}(\alpha) \in {\Bbb F}_{q^n}$.
\end{lemma}
{\bf Proof}.  Assume there is such an element $\alpha$. Then we have that 
$$ \alpha^{q^{n-2}} + \alpha^{q^{n-3}} + \ldots + \alpha \in {\Bbb F}_{q^n},$$
which gives, on raising to the power of $q$
$$ \alpha^{q^{n-1}} + \alpha^{q^{n-2}} + \ldots + \alpha^q \in {\Bbb F}_{q^n}.$$
Subtracting the first equation from the second, we get
$$ \alpha^{q^{n-1}} - \alpha \in {\Bbb F}_{q^n}. $$
This implies that 
$$ \alpha^{q^{2n-1}} - \alpha^{q^n} = \alpha^{q^{n-1}} - \alpha. $$
Now raising both sides to the power of $q^{n-1}$, we get
$$ \alpha^{q^{n-2}} - \alpha^{q^{2n-1}} = \alpha^{q^{2n-2}} - \alpha^{q^{n-1}}. $$
Adding this equation to the previous one, we get
$$ \alpha^{q^{n-2}} -  \alpha^{q^n} = \alpha^{q^{2n-2}} - \alpha. $$
This gives us 
$$ \alpha^{q^n} - \alpha \in {\Bbb F}_{q^{n-2}}. $$

For this to hold, we would need ${\Bbb F}_{q^{n-2}}$ to be a subfield of ${\Bbb F}_{q^{2n}}$, or, $n-2\;|\;2n$. But for $n-2 > 0$,  this happens only for $n=3,4\mbox{ and } 6$. \hfill $\Box$

\begin{theorem} If $n \neq 3,4 \mbox{ or } 6$, 
$$N_2(E_2) = 0.$$
\end{theorem}
{\bf Proof}. Follows from Lemma~\ref{lemma:N_2=0} and Lemma~\ref{lemma:neq-2-3-6}.  \hfill $\Box$

\vspace{0.2in}
We conclude this chapter with an interesting insight into the curves corresponding to symmetric function fields. On any curve, for rational point $(\alpha, \beta)$ on the curve, the point's conjugate points  $({\alpha} ^{q^i}, {\beta}^{q^i}), 1 \leq i \leq n-1$, also lie on the curve. However, in general,  points 
$({\alpha},{\beta}^{q^i}), 1 \leq i \leq n-1$, and $ ({\alpha}^{q^i}, {\beta}), 1 \leq i \leq n-1$, do not lie on the curve.  In the case of a curve
corresponding to a symmetric extension of the rational function field,  however, these  points too lie on the curve.  Thus we get many rational points ``for free". Similarly, we get many automorphisms also.



\chapter{Quasi-symmetry} \label{chapter:quasi-symmetry}

\section{Introduction}

In Chapter~\ref{chapter:symmetry}, we described some families of extensions of the rational function field in which all but possibly one of the places of degree one in the rational function field split completely in the extension. We showed how many existing function fields that are known to be maximal in some sense (Hasse-Weil, Oesterle) actually fall in some of these families, and we presented more, hitherto unknown families with the same amount of splitting of rational places, but lower genus. In some cases, these families contained function fields that were maximal in the Oesterle sense. All of the above was done using the notion of symmetry of functions, in the appropriate context.

In this chapter, we introduce a larger class of functions, those that we call ``quasi-symmetric,'' that can be used to the same effect. The notion of quasi-symmetry generalizes that of symmetry for purposes of producing extensions in which almost all rational places split completely. We show how, using quasi-symmetric functions, it is actually possible to split {\em all} the rational places in extensions of function fields. Thus, some of the extensions described in this chapter can attain the maximum possible value for the ratio of $N(E)/[E:F]$, for a fixed $[E:F]$. 

As is true for Chapter~\ref{chapter:symmetry}, the methods described in this chapter can be applied to construct explicit extensions of arbitrarily high genera and number of rational places.

We also demonstrate that because the set of quasi-symmetric functions is considerably larger than that of symmetric functions, it is possible to use them to control the splitting and ramification of rational places more effectively. In particular, we can use them to ``push'' all the ramification out of the set of rational places altogether.

\section{Quasi-symmetric functions}

We now introduce a notion we will use extensively. The notation will be as in 
Section~\ref{section:symmetric-functions}.

\begin{definition} A polynomial ${\bf f}(X)$ in $R[X]$ will be called quasi-symmetric if it is fixed by the cycle $\varepsilon = (1\; 2\;\ldots\; n) \in {\frak S}_n$.  If $\varepsilon$ is allowed to act on $\overline{R}(X)$ in the natural way, its fixed points will be called quasi-symmetric rational functions, or simply, quasi-symmetric functions. These form a subfield $\overline{R}(X)_{qs}$  of $\overline{R}(X)$.
\end{definition}

\begin{lemma} For $n>2$, there always  exist quasi-symmetric functions that are not symmetric. 
\end{lemma}
{\bf Proof}.  $\langle{\varepsilon}\rangle$ has index $(n-1)!$ in $\sn$. Thus for $n >2$, the set of functions fixed by \sn is strictly smaller than those fixed by $(\varepsilon)$. For $n=2,\; \sn = (\varepsilon)$ so that the notions of symmetric and quasi-symmetric coincide.\hfill $\Box$

\begin{example} \rm \label{example:quasi-symmetric-n=3}
\rm $(n=3)$ A family of quasi-symmetric functions in three variables is given below:
$${\bf f}(x_1,x_2,x_3) = x_1{x_2}^i + x_2{x_3}^i + x_3{x_1}^i. $$ 
Note that for $i \neq 0 \mbox{ or }1$, these are not symmetric.
\end{example}

\begin{definition}
Consider the  extension ${\Bbb F}_{q^n}/{\Bbb F}_q$ of finite fields cf. Proposition~$\ref{proposition:extension-finite-fields}$. We will evaluate the  quasi-symmetric polynomials (resp. quasi-symmetric functions) in ${\Bbb F}_{q^n}(X)$ at  $(X)=(t,\phi(t),\ldots,\phi^{n-1}(t))=(t,t^q,\ldots,t^{q^{n-1}})$. These will be called $(n,q)$-quasi-symmetric  polynomials (resp. $(n,q)$-quasi-symmetric functions). 
\end{definition}

\begin{example}\rm
Using the three-variable quasi-symmetric functions of Example~$\ref{example:quasi-symmetric-n=3}$, we can obtain  the following $(3,q)$-quasi-symmetric functions:
$$f(t) = {\bf f}(t,t^q,t^{q^2}) = t^{1+iq} + t^{q+iq^2} + t^{q^2 + i}.$$
Again, these are not symmetric for  $i \neq 0 \mbox{ or }1$.
\end{example}

\begin{lemma} \label{lemma:quasi-is-invariant-ifpart}
Let ${\bf f}(X) \in R[X]$. Then if ${\bf f}(X)$ is quasi-symmetric, $f(t^q) = f(t) \bmod (t^{q^n}-t)$.
\end{lemma}
{\bf Proof}.  We have that
$$ {\bf f}(\varepsilon(X)) = {\bf f}(x_2,x_3,\ldots,x_n,x_1). $$
Evaluating ${\bf f}(\varepsilon(X))$ at $(X)=(t,t^q,\ldots,t^{q^{n-1}})$, we get
$$  {\bf f}(t^q,t^{q^2},\ldots,t^{q^{n-1}},t) =  f(t^q)\bmod(t^{q^n} - t). $$
But since $f(X)$ is quasi-symmetric, this is equal to ${\bf f}(X)$ evaluated at $(X)=(t,t^q,\ldots,t^{q^{n-1}})$, which is ${\bf f}(t, t^q,\ldots,t^{q^{n-1}}) = f(t)$.   \hfill $\Box$

There is a converse to this statement if we allow only polynomials with degrees less than $q^n$. But first we need the following definition.

\begin{definition}
A lift of $f(t)$ is a polynomial ${\bf f}(X)$ such that ${\bf f}(t,t^q,\ldots,t^{q^{n-1}}) = f(t)$.
\end{definition}

\begin{lemma}  \label{lemma:quasi-is-invariant}
Let $f(t) \in \nq[t]$ have degree less than $q^n$. Then $f(t)$ is $(n,q)$-quasi-symmetric  iff $f(t^q) = f(t)\bmod (t^{q^n}-t)$.
\end{lemma}
{\bf Proof}. $\Rightarrow$ As in Lemma~\ref{lemma:quasi-is-invariant-ifpart}.

$\Leftarrow$ Each monomial in $f$ of degree $d$, with $d<q^n$ has a unique lift $x_1^{d_1}\ldots x_n^{d_n}$, where $d_i < q, 1 \leq i \leq n$ and $d_nd_{n-1}\ldots d_1$ is the base $q$ representation of $d$. 

Then the lift of $f$ is just the sum of the lifts of each of its monomials. Now define $f_q(t)$ as the unique polynomial of degree less than $q^n$ which is congruent to $f(t^q)$ modulo $(t^{q^n}-t)$. 

Now note that the lift of $f_q$ is just the cyclic shift in the variables of the lift of $f$. But by hypothesis, $f(t) = f_q(t)$ and thus the lift of $f$ must be invariant under a cyclic shift (\ie, under the action of $\varepsilon$), and thus is quasi-symmetric. It then follows that $f(t)$ is $(n,q)$-quasi-symmetric. \hfill $\Box$

\begin{corollary} \label{corollary:f-is-quasi}
A polynomial $f(t) \in \nq[t]$ is $(n,q)$-quasi-symmetric iff $f(\gamma^q) = f(\gamma), \;\; \forall \gamma \in \nq$.
\end{corollary}

\begin{lemma}
For a $(n,q)$-quasi-symmetric polynomial $f(t)$ with \fq coefficients, $f(\gamma) \in {\Bbb F}_q,\; \forall \gamma \in {\Bbb F}_{q^n}$.
\end{lemma}
{\bf Proof}. It suffices to prove that $\phi(f(\gamma))= f(\gamma)$. But this is immediate from the fact that the coefficients of $f$ are in ${\Bbb F}_q$ and  from Corollary~\ref{corollary:f-is-quasi}.  \hfill $\Box$

\begin{lemma}
Let $f$ be any function from \nq to \nq that satisfies $f(t^q) = f(t)$. Then  $f|_{\nq} = g(t)$, where $g(t)$ is a $(n,q)$-quasi-symmetric polynomial. 
\end{lemma}
{\bf Proof}. By Lagrange interpolation, we can find a polynomial $g(t)$ in $\nq[t]$ that agrees with $f$ on \nq. Thus, $g(\gamma^q) = g(\gamma), \forall \gamma \in \nq$. Now from Corollary~\ref{corollary:f-is-quasi}, it follows that $g(t)$ is $(n,q)$-quasi-symmetric.\hfill $\Box$

Due to this result, we may use the terms $(n,q)$-quasi-symmetric function and $(n,q)$-quasi-symmetric polynomial interchangeably when the context allows. 

\begin{definition}
Let $U_{qs}$  denote the ${\Bbb F}_{q^n}$-vector space of $(n,q)$-quasi-symmetric functions from \nq to $\nq\!,$ and $V_{qs} \subset U_{qs}$ denote the ${\Bbb F}_q$-space of all functions $f \in U_{qs}$ such that $f|_{\nq} =  g(t)|_{\nq}$ for some $g(t) \in {\Bbb F}_q[t]$. Thus, $V_{qs}$ consists of $(n,q)$-quasi-symmetric functions having a polynomial representation with coefficients in \fq. 
\end{definition}

\begin{lemma}
Any ${\Bbb F}_q$-linearly independent subset $\{f_i(t)\}_{1 \leq i \leq r}$ in $V_{qs}$ is also a ${\Bbb F}_{q^n}$-linearly independent subset in $V_{qs}$.
\end{lemma}
{\bf Proof}. Let $u_1f_1(t) +\dots + u_rf_r(t) = 0, u_i \in {\Bbb F}_{q^n}$. Now express the $u_i$ as $\sum a_{ij}w_j, a_{ij} \in {\Bbb F}_q, w_1,\dots,w_n$ a ${\Bbb F}_q$-basis for ${\Bbb F}_{q^n}$, then $\sum(\sum a_{ij}f_i(t))w_j = 0$ and since $f_i(u) \in {\Bbb F}_q, u \in {\Bbb F}_{q^n}$, it follows that
$\sum a_{ij}f_i(t)=0$ and hence $a_{ij} = 0$.\hfill $\Box$

\begin{definition}
Let ${\mathcal O}$ denote the set of orbits of Galois for the action of $\gal({\Bbb F}_{q^n}/{\Bbb F}_q)$ on \nq. 
\end{definition}

We note that $U_{qs}$ consists of functions from ${\Bbb F}_{q^n}$ to ${\Bbb F}_{q^n}$ which are constant on each orbit in ${\mathcal O}$. $V_{qs}$ consists of functions from \nq to \fq that are constant on these orbits.

\begin{lemma}
Let $f \in V_{qs}$, then there is a  unique polynomial representation $g(t) \in \fq[t]$, of degree $< q^n$, such that $f|_{\nq} = g(t)$.
\end{lemma}
{\bf Proof}. Choose a polynomial representation $g(t)$ of $f$ degree less that $q^n$. Define $\phi(g(t))_q$ as the unique polynomial of degree less than $q^n$ which is congruent to $\phi(g(t))$ modulo $(t^{q^n} - t)$. Since $f \in V_{qs}$,  we must have $g(t), \phi(g(t))_q$ agree on \nq. Therefore, they must be equal as their degree is less than $q^n$ (for if they were not, then the polynomial $g(t)-\phi(g(t))_q$ would have more zeros than its degree).  Thus, the coefficients of $g(t)$ are fixed by the Galois group and must lie in ${\Bbb F}_q$. Uniqueness follows similarly.\hfill $\Box$

When there is no scope for confusion, we will often say that a polynomial is in $V_{qs}$ to mean that it is the representative of a function in $V_{qs}$.

Now we examine the dimensions of $U_{qs}$ and $V_{qs}$. 

\begin{theorem}
$U_{qs}$ and $V_{qs}$ have dimension $|{\mathcal O}|$ as vector spaces over \nq and \fq respectively
\end{theorem}
{\bf Proof}.  Consider the set of functions $\{f_i\}_{1 \leq i \leq |{\mathcal O}|}$, where $f_i$ is a function that is $1$ on the \ith orbit of Galois and $0$ on all others. Then this set is linearly independent over \nq (resp. \fq), and it spans $U_{qs}$ (resp. $V_{qs}$). \hfill $\Box$

\begin{lemma} \label{lemma:qs-no-zeros}
There exist $(n,q)$-quasi-symmetric functions that have no zeros in ${\Bbb 
F}_{q^n}$.
\end{lemma}
{\bf Proof}. We are free to assign any value in \fq for a $(n,q)$-quasi-symmetric function on an orbit of Galois. In particular, we are free to assign a non-zero value to each orbit. Then there are $(q-1)^{|{\mathcal O}|}$ such functions constant on the orbits of Galois, and having no zero in \nq.\hfill $\Box$

We now discuss some ways to actually construct such functions whose existence is guaranteed by Lemma~\ref{lemma:qs-no-zeros}. We first consider a type of construction that we call the ``composition with irreducibles'' construction.

\begin{lemma} \label{lemma:composition-with-irreducibles}
Let $i(t)$ be a polynomial of degree $d_1$ over ${\Bbb F}_q$ that has no roots in \fq. In particular, we may choose $i(t)$ to be irreducible over \fq. Let $s(t) \in V_{qs}$ be a $(n,q)$-quasi-symmetric polynomial of degree $d_2$. Then the following hold:
\begin{enumerate}
\item $i(s(t))$ is $(n,q)$-quasi-symmetric, for $1 \leq i \leq n$, and maps ${\Bbb F}_{q^n}$ to ${\Bbb F}_q$.
\item $i(s(t))$ has no zeros in ${\Bbb F}_{q^n}$.
\end{enumerate}
\end{lemma}
{\bf Proof}.
For (i), note that $s(t)$  maps ${\Bbb F}_{q^n}$ to ${\Bbb F}_q$. Now since the coefficients of $i(t)$ are from \fq, we get the result.  For (ii), assume the contrary, \ie, let 
there be $\alpha \in {\Bbb F}_{q^n}$ s.t $f(s(\alpha)) = 0$. Then we have that $s(\alpha)$ is a zero of $i(t)$. But $s(\alpha) \in {\Bbb 
F}_q$, since $s(t) \in V_{qs}$. This would then give us a zero of $i(t)$ in ${\Bbb F}_q$, contradicting the hypothesis. \hfill $\Box$

The resulting polynomial obtained by this composition would have degree equal to $d_1d_2$. Since there exist irreducibles of any arbitrary degree $>1$ over \fq, we can now construct $(n,q)$-quasi-symmetric polynomials that map \nq to \fq  having degree equal to any multiple of a $(n,q)$-quasi-symmetric function in $V_{qs}$, and without a zero in ${\Bbb F}_{q^n}$.

A particularly useful form of this construction is given below.

\begin{corollary}
Let $s(t) \in V_{qs}$ be a $(n,q)$-quasi-symmetric polynomial. Then the polynomial $f(t)= {s(t)}^m - \beta$, where $\beta \in {\Bbb F}_q$ is not a $m^{th}$ power in ${\Bbb F}_q$, is $(n,q)$-quasi-symmetric, 
maps ${\Bbb F}_{q^n}$ to ${\Bbb F}_q$, and has no zeros in ${\Bbb F}_{q^n}$.
\end{corollary}
{\bf Proof}. Follows from Lemma~\ref{lemma:composition-with-irreducibles}, by choosing $i(t) = t^m - \beta$.\hfill $\Box$

\begin{example}\rm
The polynomial $f(t) = t^{10} + 2t^6 + t^2 - 2$ is $(2,5)$-quasi-symmetric and has no zeros in ${\Bbb F}_{25}$. We obtain $f(t)$ by composing the irreducible $ i(t) = t^2 - 2$ with the first $(2,5)$-elementary symmetric polynomial, given by $s_{2,1}(t) = t^5 + t$.
\end{example}

Now we wish to use these notions in the setting of the algebraic function field $F/K$, where $K =\nq$.
\begin{definition}
$F_{qs}$ will denote  the field of $(n,q)$-quasi-symmetric functions in $F$. $F^{\phi}_{qs}$ will denote the field of $(n,q)$-quasi-symmetric rational functions in $F$, whose coefficients are from \fq.
\end{definition}

\section{Quasi-symmetric extensions of function fields}

$F$ and $K$ are as described earlier. $E$ is a finite separable extension of $F$, generated by $y$, where $\varphi(y) = 0$, for $\varphi(T)$ an irreducible polynomial in $F[T]$.

In this section we will introduce families of extensions of $F$ whose generators satisfy explicit equations involving only $(n,q)$-quasi-symmetric functions. Let $y$ satisfy
                        $$ g(y) = f(x),$$
where $f, g \in F_{qs}^\phi$. If $K=\nq$, this implies that in the residue field of a rational place, alhough the class of $x$ and $y$ can assume any values in ${\Bbb F}_{q^n} \cup \infty$, that of $f(x)$ and $g(y)$ will assume values only in ${\Bbb F}_{q}\cup \infty$. Among the Galois extensions that such equations can produce are the two special cases of extensions of Artin-Schreier and Kummer type. 

\subsection{Quasi-symmetric extensions of Artin-Schreier type}

We are now in a position to state the main theorems of this section. 
\begin{theorem} \label{theorem:general-quasi-symmetric}
Let $F=K(x)$ where $K=\nq$. Let $E=F(y)$, where $y$ satisfies the equation
\begin{equation} \label{equation:general-quasi-symmetric}
      y^{q^{n-1}} + y^{q^{n-2}} + \ldots + y = \frac{h(x)}{g(x)},
\end{equation}
where $h(x),g(x) \in F_{qs}^\phi$, and $\frac{h(x)}{g(x)}$ is not the image of any rational function in $F$ under a linear polynomial. Then the following hold:
\begin{enumerate}
\item $E/F$ is a Galois extension, with degree $[E:F] = q^{n-1}$. $\gal(E/F)= \{\sigma_\beta:y \rightarrow y + \beta\}_{s_{n,1}(\beta) = 0}.$  
\item Let $P \in \pf$ be such that $v_P \left(\frac{h(x)}{g(x)}\right) = -m, \;m>0$ and $\gcd(m,q) = 1$. Then $P$ is totally ramified in $E$, with different exponent
$$ d(P'|P) = (q^{n-1} -1)(m + 1). $$
\item Let $Q \in \pf$ be any rational place, such that  $v_Q \left(\frac{h(x)}{g(x)}\right) \geq 0$. Then $Q$ splits completely in $E$. 
\end{enumerate}
\end{theorem}
{\bf Proof}. For (i), irreducibility follows from Proposition~\ref{proposition:irreducibility}. The proof of that proposition carries through even for the case of rational functions. Also, we note that since  $h(x),g(x) \in F_{qs}^\phi$, the residue class of $\frac{h(x)}{g(x)}$ at any rational place is in \fq. Now (iii) follows from  Proposition~\ref{proposition:Art-Sch}.\hfill $\Box$

\begin{example}\rm
Let $F=K(x)$ and $K = {\Bbb F}_{q^3}$. Let $E=F(y)$, where $y$ satisfies the equation  
$$ y^{q^2} + y^q + y = x^{1+iq} + x^{q +iq^2 } + x^{q^2+i}. $$
Now from Lemma~$\ref{lemma:subextensions}$, we can get a subextension $E^1$ whose degree over $F$ is $q$. Thus, $E^1 = F(y_1)$ where $y_1$ satisfies the equation
$$ y_1^q - b^{1-q}y_1 = x^{1+iq} + x^{q +iq^2 } + x^{q^2+i}, $$
where $b$ is a nonzero element of ${\Bbb F}_{q^3}$ whose trace in \fq is zero.
By substituting $y_1 = z_1 + x^{1+iq},$ we get the following equation
$$ z_1^q - b^{1-q}z_1 = (1+b^{1-q})x^{1+iq} + x^{q^2+i}. $$
In the case where $0 < i < q$, the degree of the RHS is $q^2 + i$. Also let $\gcd(i,q)=1$. Then we get that if $P_\infty^1$ is the place dividing $P_\infty$ in $E^1$,
$$ d(P_\infty^1|P_\infty) = (q -1)(q^2+i+1). $$
For the case of $i > q+1$, the degree of the RHS is $1+iq$  and then 
$$ d(P_\infty^1|P_\infty) = (q -1)(2 + iq). $$ 
In case $i=q+1$, both the terms in the RHS have equal degrees. Notice that $1 + b^{1-q} = -b^{q^2-q} \neq 1$ unless the characteristic is three. Thus, if the characteristic is not three, we have that 
$$ z_1^q - b^{1-q}z_1 = (1- b^{q^2-q})\;x^{q^2+q+1} \neq 0.$$
The different exponent is given by
$$ d(P_\infty^1|P_\infty) = (q^2 -1)(q^2+q+2). $$ 
In all these cases, we have that $N(E) = q^5+1$, since all rational places, except for $P_\infty$ split completely. For $i=1$ the RHS is the  $(3,q)$-elementary symmetric polynomial $s_{3,2}$. In that case, for $q=2$, $g(E) = 6$ and the extension attains the Oesterle lower bound on genus.
\end{example}

We now use $(n,q)$-quasi-symmetric functions that have no zeros in $K$, to build extensions of $F/K$. 
 
\begin{theorem} 
Let $F=K(x)$, where $K = \nq$. Let $E=F(y)$, where $y$ satisfies the equation 
 $$y^{q^{n-1}} + y^{q^{n-2}} + \ldots + y = \frac{h(x)}{g(x)}.$$
where $h(x), g(x) \in F_{qs}^\phi$, $\frac{h(x)}{g(x)}$ is not the image of a rational function in $F$ under a linear polynomial, $deg(g(x)) > deg(h(x))$, and $g(x)$ has no zeros in \nq. Then all the rational places of $F$ split completely in $E$. Thus $N(E) = q^{n-1}(q^n +1)$, and $N(E)/[E:F]$ attains its maximum possible value of $q^n + 1$.
\end{theorem}
{\bf Proof}. For all $\alpha \in \nq$, $ h(\alpha)/g(\alpha) \in \fq$, which ensures splitting of all rational places of the form $P_\alpha$. Also, since $deg(g(x)) > deg(h(x))$, the RHS of the equation has a zero at $P_\infty$, ensuring that $P_\infty$ also splits completely in $E$.\hfill $\Box$

\begin{example} \rm
Let $F=K(x)$, where $K = {\Bbb F}_{q^2},\; p \neq 2$. Let $E=F(y)$, where $y$ satisfies the equation
$$y^q + y = \frac{x^{q+1}}{x^{2q} + 2x^{q+1} + x^2 - \alpha},$$
where $\alpha$ is not a square in \fq. All the rational places in $F$ split completely in $E$. The non-rational places which divide $x^{2q} + 2x^{q+1} + x^2 - \alpha$ are totally ramified. 
For the special case of $q=5$, and obtaining a quadratic extension of ${\Bbb F}_5$ using the root $t$ of the irreducible polynomial $y^2 - 3$, we get the following  factorization into irreducibles of $x^{10} + 2x^{6} + x^2 - 2 $, where $2$ is our chosen non-square in  ${\Bbb F}_5$:
$$ x^{10} + 2x^{6} + x^2 - 2 = (x^5 + x + 2t)(x^5 + x + 3t). $$
This then gives us the following degree of different for $E/F$
$$ \degdiff(E/F) = 2(5-1)(1+1)(5) = 80, $$
which gives us the genus $g(E) = 37$. Also, since all rational places split completely, $N(E) = 130$. For this value of $q$ and $N$, the Oesterle lower bound on genus is $g \geq 11$. Thus while all the rational places split completely in this extension, the rise in genus that we incur in achieving this maximum splitting is high. Thus the $N/g$ ratio for $E$ is not very high. This is a typical example in this respect. Note that we can split all except one rational place of ${\Bbb F}_{25}$ in a degree $5$ extension by using the Hermitian function field, given by $E = F(y)$, where $y$ satisfies the equation
$$ y^5 + y = x^{6}.$$
In this case, the only place in \pf that is ramified is $P_\infty$, which is totally ramified with different exponent given by
$$ d(P'_\infty|P_\infty) = 5+2 = 7, \mbox{  and,} $$
$$ \degdiff(E/F) = (5-1)(6+1) = 28.$$
Thus, to split the one remaining rational place, we have had to almost triple the degree of the different. 
\end{example}

\subsection{Quasi-symmetric extensions of Kummer type}

We now study extensions whose Galois group is a subgroup of the 
multiplicative group $K^*$. For this we will need that the field 
contain a primitive \ith root of unity $\xi_i$ for some $i$ coprime to $p$. In particular we know that $K$ contains $\xi_i$ for $i = \frac{q^n-1}{q-1}$.                                 

\begin{theorem} \label{theorem:Kummer-quasi-symmetric}
Let $F=K(x)$ where $K = \nq$. Let $E=F(y)$, where $y$ satisfies the equation
               $$ y^{\frac{q^n-1}{q-1}} = \frac{h(x)}{g(x)}, $$
where  $h(x), g(x) \in  F_{qs}^\phi$ and $\frac{h(x)}{g(x)} \neq w^{\frac{q^n-1}{q-1}}, \; \forall w \in F$. Then the following hold:
\begin{enumerate}
\item This is a cyclic Galois extension, with degree $[E:F] = \frac{q^n-1}{q-1}$. $\gal(E/F) = \{ \sigma_j : y \rightarrow y{\xi}^j \}_{1 \leq j \leq \frac{q^n-1}{q-1}}$.
\item The only places of $F$ that may be ramified are the unique pole $P_{\infty}$ of $x$ and the zeros of $h(x)$ and $g(x)$. For any such place $P$, let $v_P = v_{P}(h(x)/g(x))$ be the corresponding valuation. Let $P'$ denote a generic place lying above $P$ in $E$. Then the ramification index for $P'$ over $P$ is given by 
$$ e(P'|P) =\frac{[E:F]}{r_{P}},  $$
where $r_{P} = \gcd([E:F],v_P) > 0$. Since the extension is tame, the different exponent for $P'$ over $P$ is given by
$$  d(P'|P) = e(P'|P)-1. $$ 
\item All other rational places of ${\Bbb F}_{q^n}(x)$ split completely in the 
extension. The rational places discussed in {\rm (ii)}, if not ramified, also split completely. 
\end{enumerate}
\end{theorem}
{\bf Proof}. The proofs of (i) and (ii) can be looked up in any 
reference on Kummer extensions. An excellent reference is \cite{Sti1}. (iii) follows from the fact that for $\alpha \in {\Bbb F}_{q^n}$ and  not  a zero of $h(x)$ or $g(x)$, $h(\alpha)/g(\alpha) \in {\Bbb F}_q^*$ and then we know that it has  $\frac{q^n-1}{q-1}$ pre-images under the norm map.\hfill $\Box$

\begin{corollary}
If in Theorem~$\ref{theorem:Kummer-quasi-symmetric}$, $h(x)$ and $g(x)$ have no zeros in ${\Bbb F}_{q^n}$, then other than possibly $P_{\infty}$, all rational places split completely in the extension, giving atleast $\frac{q^{2n} -q^n}{q-1}$ rational places. In addition, if $r_{\infty} = [E:F]$, then $P_{\infty}$ also splits completely. In that case, all rational places split completely in the extension.
\end{corollary}

\begin{example} \rm Let $q=p^e, e>1, p \neq 2, m= \frac{q-1}{p-1}$ and $\alpha$ not a square in ${\Bbb F}_p$. Let $F=K(x)$, where $K = \fq$, and $E=F(y)$, where $y$ satisfes the equation
 $$ y^m = x^{2m} - \alpha. $$
Now all the rational places will spilt completely. Ramification will be restricted to the non-rational places that correspond to the irreducible polynomials
which divide  $x^{2m} - \alpha$. Thus, this extension will attain the maximum possible value for the ratio $N(E)/[E:F]$, for $[E:F]=m$.
\end{example}

We have discussed a few techniques that can be used to split ``almost all'' rational places in extensions of function fields. In most of our examples, that means all, or all except one. We observe that it is very hard to keep the genus low if we split all rational places. It seems that to split the last rational place entails a fairly high increase in genus, as compared to splitting all but one rational places.



\chapter{Towers}

\section{Introduction}

Let ${\mathcal F} = (F_1,F_2,F_3,\ldots)$ be a tower of function fields, each defined over $K$. Further, we will assume that $F_1 \subseteq F_2 \subseteq F_3 \ldots$, where $F_{i+1}/F_i$ is a finite separable extension and  $g(F_i) > 1$ for some $i \geq 1$. This follows the convention of \cite{GarSti2}.

In this chapter, we apply the techniques developed in Chapters~\ref{chapter:symmetry} and \ref{chapter:quasi-symmetry} to splitting rational places in towers of function fields. While the basic ideas are the same, one has to keep in mind that  what is optimal at one stage of the tower may lead to complications at later stages. 

Let ${\mathcal F}$ be as above. It is known that the sequence $(N(F_i)/g(F_i))$ converges as $i \rightarrow \infty$ \cite{GarSti2}. Let $\lambda({\mathcal F}) := \lim_{i \rightarrow \infty} N(F_i)/g(F_i)$. 

There are known bounds on the behaviour of function fields over a finite field \fq. Let $N_q(g)$ := max\{$N(F) / F$ a function field over ${\Bbb F}_q$ of genus $g(F) = g$\}. Also, let
\begin{equation}
A(q) := \limsup_{g \rightarrow \infty} N_q(g)/g,
\end{equation}
then the Drinfeld-Vladut bound \cite{DriVla1} says that
\begin{equation}
A(q) \leq \sqrt{q} - 1.
\end{equation}
Ihara \cite{Iha1}, and Tsafasman, Vladut and Zink \cite{TsfVlaZin1} showed that this bound can be met in the case where $q$ is a square. It is not known what the value of $A(q)$ is for non-square $q$, though there are results by Serre \cite{Ser1,Ser2,Ser3} and Schoof \cite{Sch1} in this direction.

Clearly, for a tower of function fields ${\mathcal F} = (F_1,F_2,\ldots)$, $F_i/\fq$, we have that
\begin{equation}
 0 \leq \lambda({\mathcal F}) \leq A(q).  
\end{equation}

Garcia and  Stichtenoth \cite{GarSti1}, \cite{GarSti2} gave two explicitly constructed towers of function fields over a field of square cardinality that meet the Drinfeld-Vladut bound asymptotically. In \cite{GarSti3}, they gave more explicit descriptions of towers of function fields over \fq, with $\lambda({\mathcal F}) > 0$. These also meet the Drinfeld-Vladut bound in some cases where the underlying field of constants is of square cardinality.

Elkies, in \cite{Elk1}, gave eight explicit iterated equations for towers of modular curves, which also attained the Drinfeld-Vladut bound over certain fields and showed that the examples presented in \cite{GarSti1} and \cite{GarSti3} were also modular. He then conjectured that all asymptotically optimal towers would, similarly, be modular. 

In Chapter~\ref{chapter:symmetry}, we used the notion of symmetry of functions to describe explicitly constructed extensions of function field in which all rational places except one split completely. In Chapter~\ref{chapter:quasi-symmetry}, we showed that on generalizing the notion of symmetry to include the so-called ``quasi-symmetric'' functions, we could actually split all the rational places in an extension of function fields. Furthermore, in both these cases, we obtained infinite families of extensions with such properties.

In this chapter, we apply the techniques developed in Chapter~\ref{chapter:symmetry} and Chapter~\ref{chapter:quasi-symmetry} to the problem of splitting rational places in a tower of function fields. Towards that end, we describe infinite families of towers in which all the rational places split completely throughout the tower. We also describe infinite families of towers in which all rational places, except one, split completely throughout the tower. Then we observe that inspite of such splitting behaviour at the rational places, all these towers have $\lambda({\mathcal F}) = 0$. In that sense, the main accent here is not so much on obtaining a high value for  $\lambda({\mathcal F})$, as it is to show the existence of certain explicitly constructed families of towers in which all rational places split completely throughout the tower. In addition, it is hoped that these examples will lead to a better general understanding of what makes  $\lambda({\mathcal F}) > 0$.  We also generalize two examples of towers with  $\lambda({\mathcal F}) > 0$ presented in \cite{GarSti3} to obtain infinite families of such towers. Subfamilies of these attain the Drinfeld-Vladut bound.

\section{Families of towers attaining the Drinfeld-Vladut bound}

\begin{theorem}  \label{theorem:good-family-1}
Let $q=p^n$ and $m|n, m \neq n$. Let $k_m = (p^n - 1)/(p^m -1)$.  Consider a tower of function fields in the family given by ${\mathcal T} = (T_1, T_2, \ldots)$, where $T_1 = \fq(x_1)$ and for $i \geq 1$, $T_{i+1} = T_i(x_{i+1})$, where $x_{i+1}$ satisfies
\begin{eqnarray*}
 x_{i+1}^{k_m} + z_i^{k_m} = b_i^{k_m}, \\
                          z_i = a_ix_i^{r_i} + b_i,
\end{eqnarray*}
where $a_i,b_i \in {\Bbb F}_{p^m} \setminus \{0\}$ for $i \geq 1$. Also $r_i$ is a power of $p, \;\forall i$. Then the following hold:
\begin{enumerate}
\item $P_\infty$ splits completely throughout the tower.
\item Every ramified place in the tower lies above a rational place in $T_1$.
\item $\lambda({\mathcal T}) \geq \frac{2}{q-2}$, and hence this family attains the Drinfeld-Vladut bound for $n = 2, m=1$ and $q=4$.  
\end{enumerate}
\end{theorem}
{\bf Proof}. 
Firstly, we verify that under the hypothesis, we do indeed get a tower of function fields.   Notice that at one of the places dividing $x_1$ in $T_2$, we get a zero of $x_2$ of order not divisible by $k_m$. This implies that the RHS, for $i=1$, is not of the form $w^{k_m}$, for $w \in T_1$. Further, one of the places dividing $x_2$ in $T_3$ also exhibits the same performance, and so on up the tower. Thus, each equation is irreducible and gives us an extension. 

(i) follows from the basic theory of Kummer extensions cf. \cite{Sti1}, Ch. III.7. It is important to note that linear transformations fix the place at infinity, so that it splits at each stage of the tower. 

For (ii), working with residue classes, note that for ramification to take place at the $i^{th}$ step of the tower, the norm of $z_i$ should be an element of ${\Bbb F}_{p^m}$. Thus $z_i$ must be in \fq.  Since $z_i$ is obtained by a linear tranformation with \fq coefficients of a characteristic power of $x_i$, it follows that $x_i$ must be in \fq. But the relations between the variables $x_i$ and $z_{i-1}$ at the previous step of the tower then force $z_{i-1}$, and therefore $x_{i-1}$ to be in \fq. Proceeding this way to the first step of the tower, we get that $x_1 \in \fq$. Thus every ramified place in $T_i$ divides a rational place ($\neq P_\infty$) in $T_1$. 

To get (iii), notice that 
$$ N(T_j) > k_m^j, \mbox{ for } j \geq 1. $$
Also, the degree of the different at the \jth stage of the tower is always less than the value it would have had all $q$ finite rational places ramified from the second stage of the tower onwards. Now, using the transitivity of the different, we can say that 
\begin{eqnarray*}
\degdiff(T_j/T_1) &<& q(k_m -1)[1 + k_m + \ldots + k_m^{j-2}], \\
		  & < & q (k_m^{j-1} -1).
\end{eqnarray*} 
Now using the Hurwitz-genus formula, it follows that
$$ g(T_j) < \frac{(q-2)(k_m^{j-1} -1)}{2}. $$
Giving us 
$$\lim_{j \rightarrow \infty}  N(T_j)/g(T_j) \geq \frac{2}{q-2}. $$
\hfill $\Box$

This tower, for the case of $m=r_i=1; \, z_i = x_i+1$, first appeared in \cite{GarSti3}. 

\begin{theorem} \label{theorem:good-family-2}
Let $q=p^n > 4$ and $m|n$. Let $l_m = (p^m -1)$. Consider a tower of function fields in the family given by ${\mathcal T} = (T_1, T_2, \ldots)$, where $T_1 = \fq(x_1)$ and for $i \geq 1$, $T_{i+1} = T_i(x_{i+1})$, where $x_{i+1}$ satisfies
\begin{eqnarray*}
x_{i+1}^{l_m} + z_i^{l_m} = 1,  \\
			z_i  = a_ix_i^{s_i} + b_i,
\end{eqnarray*}
where $a_i,b_i \in {\Bbb F}_{p^m} \setminus \{0\}$ for $i \geq 1$. Also $s_i$ is a power of $p$, $\forall i$. Then the following hold:
\begin{enumerate}
\item $P_\infty$ splits completely throughout the tower. 
\item Every ramified place in the tower lies above a rational place in $T_1$ of the form $P_\gamma$, with $\gamma \in {\Bbb F}_{q^m}$.
\item $\lambda({\mathcal T}) \geq \frac{2}{l_m - 1}$, and hence every member of this family attains the Drinfeld-Vladut bound for $n=2, m=1$ and $q=9$.  
\end{enumerate}
\end{theorem} 
{\bf Proof}. First we verify as in the proof of Theorem~\ref{theorem:good-family-1} that we do indeed get a tower of function fields. For this, note that $b_i^{l_m} = 1$. Again (i) follows from the basic theory of Kummer extensions. For (ii), we note that to have ramification at the $i^{th}$ stage of the tower, we must have that $z_i^{l_m} = 1$ implying that $z_i \in {\Bbb F}_{q^m} \setminus \{0\}$. Then by similar reasoning as in the proof of Theorem~\ref{theorem:good-family-1} , it follows that such a ramified place would divide a rational place in $T_1$ of the form $P_\gamma$, with $\gamma \in {\Bbb F}_{q^m}$. Using the Hurwitz genus formula and the transitivity property of the different along similar lines as in the proof of Theorem~\ref{theorem:good-family-1}, we get (iii). \hfill $\Box$

This tower, for the case of $s_i=m=1$, also first appeared in \cite{GarSti3}. 

Following the conjecture of Elkies, it is very likely that many of  these towers are modular. In that case, there seems to be a definite relation between some modular towers and certain symmetric towers (the other optimal constructions from \cite{GarSti1} and \cite{GarSti2} are also symmetric, and are modular as shown in \cite{Elk1}). An interesting study would be to understand under what conditions can a modular tower be written down in terms of symmetric equations.

\section{Towers where almost all rational places split completely}

In this section we construct families of towers of function fields with very good splitting behaviour. In some of the families, all rational places split completely throughout the tower, and in others, all rational places, except one, split completely throughout the tower.

\subsection{Towers of Artin-Schreier extensions}

First we begin with a tower of function fields in which all rational places split completely throughout the tower. Following the notation of Chapter~\ref{chapter:quasi-symmetry}, we will denote the subfield of $F_i$ comprising $(n,q)$-quasi-symmetric functions in $x_j$ by $F_{j,qs}$ and the subfield comprising the $(n,q)$-quasi-symmetric functions of $x_j$ with \fq coefficients by $F_{j,qs}^\phi$. In particular, in $F_i$, $F_{i,qs}^\phi$ will denote the subfield of $(n,q)$-quasi-symmetric functions of $x_i$ with \fq coefficients.

\begin{theorem} \label{theorem:all-rational-split-1}
Consider the tower of function fields ${\mathcal F} = (F_1,F_2,\ldots)$ where $F_1 = {\Bbb F}_{q^n}(x_1)$ and for $i \geq 1$, $F_{i+1} = F_{i}(x_{i+1})$, where $x_{i+1}$ satisfies the equation
\begin{equation}
 x_{i+1}^{q^{n-1}} +  x_{i+1}^{q^{n-2}} + \ldots + x_{i+1} = \frac{g(x_i)}{h(x_i)},
\end{equation}
where $g(x_i), h(x_i) \in F_{i,qs}^\phi$, $ \frac{g(x_i)}{h(x_i)}$ is not the image of a rational function under a linear polynomial, and $h(x)$ has no zeros in \nq. Also, $deg(g(x_i)) \leq deg(h(x_i))$. Then the following hold:
\begin{enumerate}
\item All the rational places of $F_1$ split completely in all steps of the tower.
\item For every place $P$ in $T_i$ that is ramified in $T_{i+1}$, the place $P'$  in $T_{i+1}$ that lies above $P$ is unramified in $T_{i+2}$. Thus, ramification at a place cannot ``continue'' up the tower. 
\end{enumerate}
\end{theorem}
{\bf Proof}. $P_\infty$ splits completely because of the condition of the degrees of $g$ and $h$. Also, the RHS is in the valuation ring at every rational place since $h$ has no zeros in \nq and $deg(g) \leq deg(h)$.  Also, its class in the residue class field is in \fq at each of these places, since the RHS is in $F_{i,qs}^\phi$. Then Proposition~\ref{proposition:Art-Sch} tells us that every rational place in $F_i$ splits completely in $F_{i+1}$. For (ii), note that if $P \in \pti$ is ramified in $T_{i+1}$, and $P'$ is a place lying above it in $T_{i+1}$, then the RHS of the equation for $x_{i+2}$ has a zero at $P'$, because of the condition on the degrees of $h$ and $g$. Thus $P'$ will be unramified in $T_{i+2}$. \hfill $\Box$

\begin{example} \rm
Consider the tower of function fields ${\mathcal F} = (F_1,F_2,\ldots)$ where $F_1 = {\Bbb F}_{q^3}(x_1)$, $q$ is  not a power of $2$, and for $i \geq 1$, $F_{i+1} = F_{i}(x_{i+1})$, where $x_{i+1}$ satisfies the equation
\begin{equation}
x_{i+1}^{q^2} + x_{i+1}^q + x_{i+1} = \frac{x_i^{2q^2 + 2q + 2}}{(x_i^{q^2} + x_i^q + x_i)^2 - \alpha_i},
\end{equation}
where $\alpha_i \in \fq$ is not a square.  All rational places split completely throughout the tower. Let the place $P$ of $T_1$ be a simple pole of the RHS in $T_1$ (\ie, for the case $i=1$.). Then, the place $P^{(i)}$ of $T_i$, where $i \geq 2$, which divides $P$, is a pole of $x_i$ of order $2^{i-2} \bmod q$. Also  notice that there will always exist such places, if we look at the equation over \clo. Thus the equation is absolutely irreducible at each stage.
\end{example}

\begin{theorem} \label{theorem:abelian-tower}
Consider the tower of function fields ${\mathcal F} = (F_1,F_2,\ldots)$ where $F_1 = {\Bbb F}_{q^n}(x_1)$, $p\neq 2$, and for $i \geq 1$, $F_{i+1} = F_{i}(x_{i+1})$, where $x_{i+1}$ satisfies the equation
\begin{equation}
 x_{i+1}^{q^{n-1}} +  x_{i+1}^{q^{n-2}} + \ldots + x_{i+1} = \frac{1}{(x_{i}^{q^{n-1}} +  x_{i}^{q^{n-2}} + \ldots + x_{i})^2 - \alpha},
\end{equation}
where $\alpha \in \fq$ is not a square. Then the following hold:
\begin{enumerate}
\item $T_i/T_1$ is an Abelian extension for $i \geq 2$.
\item All rational places split completely throughout this tower.
\item When a (non-rational) place $P \in \pti$ is ramified in $T_{i+1}$, from then on, it behaves like a rational place for splitting, and therefore splits completely further throughout the tower.
\end{enumerate}
\end{theorem}
{\bf Proof}. First we note that the equations defining the tower at each stage are indeed irreducible. For this, note that if $P$ is a place in $T_i$ that is a zero of $ (x_{i}^{q^{n-1}} +  x_{i}^{q^{n-2}} + \ldots + x_{i})^2 - \alpha)$ in $T_i$, the zero can be of degree at most two. This can be seen as follows. Let  $\sqrt{\alpha}$ be one of the square roots of $\alpha$. Then,
\begin{eqnarray*}
x_{i}^{q^{n-1}} +  x_{i}^{q^{n-2}} + \ldots + x_{i} - \sqrt{\alpha} & = &  \frac{1}{(x_{i-1}^{q^{n-1}} +  x_{i-1}^{q^{n-2}} + \ldots + x_{i-1})^2 - \alpha} - \sqrt{\alpha}, \\
& = & \frac{1 - \sqrt{\alpha}((x_{i-1}^{q^{n-1}} +  x_{i-1}^{q^{n-2}} + \ldots + x_{i-1})^2 - \alpha)}{(x_{i-1}^{q^{n-1}} +  x_{i-1}^{q^{n-2}} + \ldots + x_{i-1})^2 - \alpha}.
\end{eqnarray*}
Now note that the second derivative of the numerator of the RHS  with respect to $x_{i-1}$ is constant. The denominator is a unit at this place. Thus the zeros of the RHS can occur to at most multiplicity two. Since a similar argument holds at each stage, the valuation of the RHS at $P$ must be a power of two, which is coprime to the characteristic.  Irreducibility then follows from Proposition~\ref{proposition:Art-Sch}. 
For (i), notice that the automorphisms of $T_{i+1}/T_i$ in the tower leave $x_{i+2}$ fixed, for $i \geq 1$. Further, $T_{i+1}/T_i$ is Abelian. For (ii), note that the class of the RHS in the residue field at any rational place is in \fq at any stage of the tower. And thus the defining equation splits into linear factors over the residue class field.  \hfill $\Box$ 

\begin{theorem} \label{theorem:wildly-ramified-all-split}
There exist wildly ramified extensions of the rational function field over non-prime fields of cardinality $> 4$ of degree equal to any power of the characteristic in which all the rational places split completely. 
\end{theorem}
{\bf Proof}. For finite-separable extensions, which are not necessarily Galois, refer Theorem~\ref{theorem:all-rational-split-1}. Each extension $T_{i+1}/T_i$ has subextensions of degree equal to any arbitrary power of $p$. By an appropriate resolution of the tower, we can get the desired result.

\begin{theorem}
There exist Abelian extensions  over non-prime fields of odd characteristic of degree equal to any power of the characteristic in which all the rational places split completely. 
\end{theorem}
For Abelian extensions,  Theorem~\ref{theorem:abelian-tower} says that the Galois group of the extension $T_i/T_1$ is an elementary Abelian group of exponent $p$, for $i \geq 1$. Thus, it will have normal subgroups of all indices that are powers of $p$. The result then follows by considering the fixed fields of these subgroups.

\begin{example} \rm
Consider the tower of function fields ${\mathcal F} = (F_1,F_2,\ldots)$ where $F_1 = {\Bbb F}_{q^3}(x_1)$, $q$ is  not a power of $2$, and for $i \geq 1$, $F_{i+1} = F_{i}(x_{i+1})$, where $x_{i+1}$ satisfies the equation
\begin{equation}
x_{i+1}^{q^2} + x_{i+1}^q + x_{i+1} = \frac{1}{(x_i^{q^2} + x_i^q + x_i)^2 - \alpha},
\end{equation}
where $\alpha \in \fq$ is not a square. In this example, all rational places split completely at all steps of the tower. Furthermore, when a (non-rational) place $P \in \pti$ is ramified in $T_{i+1}$, from then on, it behaves like a rational place for splitting, and therefore splits completely further throughout the tower.
\end{example}

\begin{theorem}
Consider the tower of function fields ${\mathcal F} = (F_1,F_2,\ldots)$ where $F_1 = {\Bbb F}_{q^n}(x_1)$ and for $i \geq 1$, $F_{i+1} = F_{i}(x_{i+1})$, where $x_{i+1}$ satisfies the equation
\begin{equation}
 x_{i+1}^{q^{n-1}} +  x_{i+1}^{q^{n-2}} + \ldots + x_{i+1} = \frac{g(x_i)}{h(x_i)},
\end{equation}
where $g(x_i), h(x_i) \in F_{i,qs}^\phi$, $ \frac{g(x_i)}{h(x_i)}$ is not linear, and $h(x)$ has no zeros in \nq. Also, $deg(g(x_i)) > deg(h(x_i))$. Then all the rational places of $F_1$, except $P_\infty$, split completely in all steps of the tower.

If, in addition, we have that  $deg(g(x_i)) = deg(h(x_i)) + 1$, the the pole order of $x_i$ in the unique place lying above $P_\infty$ in $T_i$ remains one for all $i \geq 1$.
\end{theorem}

\begin{example}\rm
Consider the tower of function fields ${\mathcal F} = (F_1,F_2,\ldots)$ where $F_1 = {\Bbb F}_{q^3}(x_1)$, $q$ is  not a power of $2$, and for $i \geq 1$, $F_{i+1} = F_{i}(x_{i+1})$, where $x_{i+1}$ satisfies the equation
\begin{equation}
  x_{i+1}^{q^2} + x_{i+1}^q + x_{i+1} = \frac{x_i^{2q^2 + 1} + x_i^{2+q} + 
x_i^{2q + q^2}}{(x_i^{q^2} + x_i^q + x_i)^2 - \alpha},
\end{equation}
where $\alpha \in \fq$ is not a square.  Here, except the unique pole $P_\infty$ of $x_1$  in $F_1$, all other rational places split completely throughout the tower. Furthermore, let $P$ be any pole of $x_2$ in $T_2$, and $P^{(n)}$ denote the unique place in $T_n$ lying above it.  Then, the  pole order of $x_n$ at $P^{(n)}$ remains constant for $n \geq 2$. 
\end{example}

\begin{example} \rm 
Consider the tower of function fields ${\mathcal F} = (F_1,F_2,\ldots)$ where $F_1 = {\Bbb F}_{q^n}(x_1)$ that is obtained as follows. $T_2 = T_1(x_1)$, where
$$ x_{2}^{q^{n-1}} +  x_{2}^{q^{n-2}} + \ldots + x_{2} = \frac{1}{(x_{1}^{q^{n-1}} +  x_{1}^{q^{n-2}} + \ldots + x_{1})^m - \alpha},
 $$
where $\alpha$ is not an $m^{th}$ power in \fq. And for $i \geq 2$, $T_{i+1} = T_i(x_{i+1})$ where $x_{i+1}$ satisfies the equation
$$ x_{i+1}^{q^{n-1}} +  x_{i+1}^{q^{n-2}} + \ldots + x_{i+1} = \frac{h(x_i)}{g(x_i)},$$
where $h(x_i),g(x_i) \in F_{i,qs}^\phi$, and $deg(h(x_i))=deg(g(x_i))+1$. Note that we are guaranteed the existence of such polynomials $h$ and $g$ by the following construction. Take any two functions $f_1$ and $f_2$ in $F_{i,qs}^\phi$ with coprime degrees $d_1$ and $d_2$ respectively (in particular, trace and norm will do). Then there exist
integers $m,n$ such that $md_1 + nd_2 = 1$. Without loss of generality, let $m$ be positive and $n$ negative. Then let $h(x_i) = i_1(f_1(x_i))$ and $g(x_i) = i_2(f_2(x_i))$, where $i_1$ and $i_2$ are irreducible polynomials over \fq of degrees $m$ and $n$ respectively. Let $P$ be any place in $T_1$ such that  $v_P({(x_{1}^{q^{n-1}} +  x_{1}^{q^{n-2}} + \ldots + x_{1})^m - \alpha}) = 1$. Then $P^{(i)}$, which is the unique place in $T_i$ dividing $P$, remains a simple pole of $x_i$ for $ i \geq 2$, ensuring irreducibility of the defining equation at each stage of the tower. 
\end{example}
 
\subsection{Towers of Kummer extensions}

\begin{theorem} \label{theorem:tamely-ramified-all-split}
Consider the tower of function fields ${\mathcal F} = (F_1,F_2,\ldots)$ where $F_1 = {\Bbb F}_{q^n}(x_1)$ and for $i \geq 1$, $F_{i+1} = F_{i}(x_{i+1})$, where $x_{i+1}$ satisfies the equation
\begin{equation}
 x_{i+1}^{\frac{q^n -1}{q-1}} = \frac{g(x_i)}{h(x_i)},
\end{equation}
where $g(x_i), h(x_i) \in F_{i,qs}^\phi$, $ \frac{g(x_i)}{h(x_i)} \neq w^{\frac{q^n -1}{q-1}}, \; \forall w \in F_i$, and $g,h$ have no zeros in \nq. Also, $deg(g(x_i)) = deg(h(x_i))$. Then all the rational places split throughout the tower.
\end{theorem}
{\bf Proof}. The RHS is in the valuation ring at every rational place at $F_i,\; \forall i$ and its class in the residue class field is in $\fq \setminus \{0\}$, since $g,h$ have no zeros in \nq, and the RHS is $(n,q)$-quasi-symmetric. Then every rational place in $F_i$ splits completely in $F_{i+1}, \; \forall i$. \hfill $\Box$

\begin{example}\rm
Consider the tower of function fields ${\mathcal F} = (F_1,F_2,\ldots)$ where $F_1 = {\Bbb F}_{q^3}(x_1)$, $q$ is  not a power of $2$, and for $i \geq 1$, $F_{i+1} = F_{i}(x_{i+1})$, where $x_{i+1}$ satisfies the equation
\begin{equation}
  x_{i+1}^{q^2 + q + 1} = \frac{(x_i^{q^2} + x_i^q + x_i)^2 - \beta}{(x_i^{q^2} + x_i^q + x_i)^2 - \alpha},
\end{equation}
where $\alpha, \beta \in \fq$ not squares. All rational places split completely throughout the tower. 
\end{example}

\begin{theorem}
Consider the tower of function fields ${\mathcal F} = (F_1,F_2,\ldots)$ where $F_1 = {\Bbb F}_{q^n}(x_1)$ and for $i \geq 1$, $F_{i+1} = F_{i}(x_{i+1})$, where $x_{i+1}$ satisfies the equation
\begin{equation}
 x_{i+1}^{\frac{q^n -1}{q-1}} = \frac{g(x_i)}{h(x_i)},
\end{equation}
where $g(x_i), h(x_i) \in $ are two $(n,q)$-quasi-symmetric polynomials, $ \frac{g(x_i)}{h(x_i)} \neq w^{\frac{q^n -1}{q-1}}, \forall w \in F_i$, and $g,h$ have no zeros in \nq. Also, $deg(g(x_i)) \neq deg(h(x_i))$. Then all the rational places, except possibly $P_\infty$ split throughout the tower.
\end{theorem}
{\bf Proof}. The RHS is in the valuation ring at every rational place in $F_i \; \forall i$, except possibly those dividing $P_\infty \in T_1$  and its class in the residue class field is in $\fq \setminus \{0\}$, since $g,h$ have no zeros in \nq, and the RHS is $(n,q)$-quasi-symmetric. Then every rational place in $F_i$ splits completely in $F_{i+1}, \; \forall i$. \hfill $\Box$

\begin{theorem} There exist tamely ramified extensions of the rational function field over a non-prime field of cardinality $>4$ of arbitrarily high degree, in which all the rational places split completely.
\end{theorem}
{\bf Proof}. Consider Theorem~\ref{theorem:tamely-ramified-all-split}. We can guarantee that such $(n,q)$-quasi-symmetric functions exist, for $q>2$. Then, in the tower described in the theorem, one can go up the tower to get arbitrarily high degree extensions of the rational function field. These will not be Galois, in general.

For the towers ${\mathcal F}$ described in this chapter in which all, or all except one, rational places split completely throughout the tower, 
$\lambda({\mathcal F}) = 0$. This is because while the ramification in the rational places is nil, or minimal, that in the non-rational places rises quite fast, leading to a fast rise in the genus. Indeed, it seems from the known examples of towers ${\mathcal F}$ with $\lambda({\mathcal F}) > 0$ that it might be necessary to have a certain amount of ramification in the rational places, in order to have $\lambda({\mathcal F}) > 0$. Or at least it seems that it is not easy to control ramification in the non-rational places, and so it is better to restrict it to a few rational places alone\footnote{These statements are for towers whose first stage is a function field of genus zero.}.

\noindent{{\bf Note:} In most of the examples that appear in this chapter, we have composed the trace/norm polynomials with the irreducible polynomial $x^2 - \alpha$, where $\alpha \in \fq$ is not a square. However, we could get infinite families of further examples by using the composition $i(q(x))$, where $i(x) \in \fq[x]$ has no zeros in \fq, and $q(x)$ is a $(n,q)$-quasi-symmetric function with \fq coefficients.



\chapter{Basis}

\section{Introduction}

The Gilbert-Varshamov bound is commonly used as a yardstick to measure the 
performance of long codes. By this measure, the commonly used BCH codes are poor at large lengths. The length of a Reed-Solomon code is limited by the size  
of its alphabet; thus large lengths require very large alphabets,
thereby greatly increasing the complexity of encoding and decoding these
codes.  It is known that there exist long alternant and concatenated codes
that meet the Gilbert-Varshamov  bound.  However, no explicit
description of these codes exists. 

Around 1980, Goppa \cite{Gop1} showed how to construct codes on algebraic curves. If $F$ is the function field corresponding to this curve, code performance depended upon the ratio $N(F)/g(F)$.  Good codes result in cases where the ratio $N(F)/g(F)$ is large and the Drinfeld-Vladut bound places an upper bound on the value of this ratio.

In 1982, using techniques from modular curves, Tsfasman, Vladut and Zink \cite{TsfVlaZin1} showed the existence of curves whose $N/g$ ratio achieved the Drinfeld-Vladut bound.  It turned out that codes built on these curves would have a performance exceeding the Gilbert-Varshamov bound -- a feat that until then was widely considered unattainable. 

However, the Tsfasman-Vladut-Zink result is existential in nature and
does not provide an explicit description of the curves. 
In 1996, Garcia and Stichtenoth \cite{GarSti1}, \cite{GarSti2} succeeded in showing that two families of curves having an explicit and simple description, also achieve the Drinfeld-Vladut bound.  However, an explicit description
of codes constructed on these curves requires the
determination of a basis for the vector spaces ${\mathcal L}(rP_\infty)$, which 
comprise functions having poles only at infinity, with pole order there bounded by $r$.

In this Chapter, we describe partial results obtained while attempting to 
construct a basis for these spaces.

\section{The discriminant, different and integral bases}

For $E/F$ a finite separable extension of function fields, let $\{\sigma_i\}_{1 \leq i \leq n}$ be the distinct embeddings of $E$ in an algebraic closure of $F$, with $\sigma_1$ being the identity. Then we have that the norm and trace of an element (resp. divisor) in $E$, is an element (resp. divisor) in $F$ given by
$$ N_{E/F}(\alpha) = \prod_{i=1}^m \sigma_i(\alpha) \mbox{ and } Tr_{E/F}(\alpha) = \sum_{i=1}^m \sigma_i(\alpha),$$
respectively.

\begin{definition}  
The field discriminant\footnote{Sometimes also known as the {\em 
discriminant of the extension.}} of $E/F$ is defined as $N_{E/F}(\diff(E/F))$, and it is a divisor in $F$. It is denoted by $\disc(E/F)$.
\end{definition}

\begin{definition}  
The set discriminant for a set of elements 
$(\alpha_1,\ldots,\alpha_n)$, with $\alpha_i \in E$, denoted $\disc(\alpha_1,\ldots,\alpha_n)$, is defined as 
\[ 
\disc(\alpha_1,\ldots,\alpha_n)=det^2\left[ \begin{array}{rrrr}
               \alpha_1 & \ldots & \alpha_n \\
      \sigma_2(\alpha_1) & \ldots & \sigma_2(\alpha_n) \\
               \vdots & \vdots & \vdots \\
      \sigma_n(\alpha_1) & \ldots & \sigma_n(\alpha_n) \\
               \end{array} \right].  \]
\end{definition}

\begin{definition}  
Let  $P \in \pf$. Let the integral closure of \op in $E$ be denoted by \oep. Then, \oep is an \op-module and there exists a basis $\{u_1,\dots,u_n\}$ of $E/F$ such that 
$$ \oep = \sum_{i=1}^n\op.u_i .$$
Any such basis will be called a local integral basis at $P$.
\end{definition}

In all the following, let $S$ be any proper subset of \pf.
\begin{definition}  We let
\begin{eqnarray*}
           & \of = \{ x \in F: v_P(x) \geq 0, \forall P \in S\},  \\
\mbox{and} & \oe = \{ x \in E: v_{P'}(x) \geq 0, \forall P' \in \pe, P' \supseteq P \in S 
\}.
\end{eqnarray*}
Then, \of (resp. \oe) will be called the ring of regular functions in $F$ 
(resp. $E$), w.r.t S. Furthermore, \oe is the integral closure of \of in $E$.
\end{definition}

\begin{definition}  
If \of is a principal ideal domain,  then there exists a basis 
$\{u_1,\ldots,u_n\}$ for \oe over \of. Thus, 
$$ \oe = \sum_{i=1}^n\of.u_i. $$
Such a basis will be called a global integral basis for \oe over \of.
\end{definition}
 
Now we recall some known facts about local and global integral bases \cite{AleDeoVijSti1}, \cite{Jan1}, \cite{Wei1}. 

\begin{proposition} \label{proposition:global=local-everywhere}
$\{\alpha_1,\ldots,\alpha_n\}$ is an integral basis for \oe 
over \of iff it is an integral basis for \oep over \op, $\forall P \in \pf$.
\end{proposition}

\begin{proposition}  \label{proposition:identify-local}
Let $\{\alpha_1,\ldots,\alpha_n\}$ be any set of elements 
in $E$ that are integral at $P$. Let $\{\beta_1,\ldots,\beta_n\}$ be a known local integral basis for \oep over \op. Then the following are equivalent:
\begin{enumerate}
\item $\{\alpha_1,\ldots,\alpha_n\}$ also form an integral basis at $P$.
\item $\disc(\alpha_1,\ldots,\alpha_n) = 
a^2\disc(\beta_1,\ldots,\beta_n)$, where $a$ is a unit at $P$.
\item $v_P(\disc(\alpha_1,\ldots,\alpha_n)) = 
v_P(\disc(\beta_1,\ldots,\beta_n)).$
\end{enumerate}
\end{proposition}

Now we state a known result that we use to identify a local integral basis at each place.

\begin{theorem} \label{theorem:setdisc=fielddisc}
Let $v_P(\beta_i) \geq 0,\; 1 \leq i \leq n$. Then  
$\{\beta_1,\ldots,\beta_n\}$ is 
a local integral basis at $P$ iff $v_P(\disc(\beta_1,\ldots,\beta_n)) = 
v_P(\disc(E/F))$.
\end{theorem}

This theorem says that a set of elements integral at $P$ form a local integral 
basis at $P$ iff their discriminant is equal to the field discriminant locally 
at $P$. 

\section{Computing the field discriminant} 

From now on, we will use the preceeding theory to search for a global integral basis for the tower meeting the Drinfeld-Vladut bound described in \cite{GarSti2}. We recall the definition of the tower below.

The tower of function fields ${\cal T} = (T_1,T_2,\ldots)$ is given by $T_1= \sq(x_1)$, and for $i \geq 2$, $T_i = T_{i-1}(x_i)$, where $x_i$ satisfies the equation 
$$ x_i^q + x_i \ = \ \frac{x_{i-1}^{q}}{x_{i-1}^{q-1}+1}. $$

Let the function $x_i^{q-1}+1$ be called $g_i$. The ramified places in $T_n$  lie above either $P_\infty$, $P_\alpha$ (where $\alpha$ is a 
zero of $g_1$) or $P_0$. The unique place in $T_n$ lying above $P_\infty$ (resp. $P_\alpha$), is denoted $P_\infty^{(n)}$ (resp. $P_\alpha^{(n)}$). The places lying above $P_0$ in $T_n$ are divided into $n$ sets $Q_i^n$, $1\leq i \leq n$, as follows. Let the unique place in $T_i$ that is a common zero of $x_1,x_2,\ldots, x_i$ be denoted by $Q_i$. Then $Q_i^n$ is the set of all places in $T_n$ that lie above $Q_i$.

We now compute the field discriminant locally at $P_\infty$, $P_\alpha$ and  $P_0$, and the patch together to get the (global) field discriminant. In each of these, we use the fact that the discriminant is the norm of the different\footnote{Pl. refer to \cite{Cas1} or \cite{Ser1}}, which, in turn, is known to us. 
\begin{enumerate}
\item The place $P_\infty$ is totally ramified throughout the tower. Thus the 
different exponent $d(P_\infty^{(n)}|P_\infty)$ can be computed from 
the transitivity of the different formula as 
$$ d(P_\infty^{(n)}|P_\infty) = 2(q-1)(1 + q + \ldots + q^{n-2}) = 2(q^{n-1} - 
1). $$
Therefore the discriminant localized at $P_\infty$ is given by  
$P_\infty^{2(q^{n-1} -1)}$.
\item The $q-1$ places $\{P_\alpha\}_{g_1(\alpha)=0}$ are also totally 
ramified all the way up the tower. Thus, from a similar computation as above, we can compute the discriminant localized at each $P_\alpha$ as being $P_\alpha^{2(q^{n-1}-1)}$. 
\item  Above $Q_k$, there are $(q-1)q^{k}$ places that ramify from $T_{2k+1}$ 
onwards.  The different exponent due to these places is, therefore, 
$2q^k(q-1)(q^{n-2k-1} -1)$. To get the discriminant localized at 
the place $P_0$, we sum these exponents over $k$ and get
\begin{eqnarray*}
\mbox{ Total Exponent} & = & \sum_{k=1}^{\lfloor\frac{n-2}{2}\rfloor} 
\mbox{Exponent 
for $Q_k$}, \\
& = &  \sum_{k=1}^{\lfloor\frac{n-2}{2}\rfloor} 2q^k(q-1)(q^{n-2k-1}-1), \\
& = & 2(q^{n-1} - q^{n- \lfloor\frac{n}{2}\rfloor} - q^{\lfloor \frac{n}{2} 
\rfloor } + q), \\
\mbox{(for $n$ even)} & = & 2(q^{n-1} - 2q^{\frac{n}{2}} + q).
\end{eqnarray*}
Thus, the discriminant localized at $P_0$ is $P_0^{2(q^{n-1} - q^{n- \lfloor\frac{n}{2}\rfloor} - q^{\lfloor \frac{n}{2} 
\rfloor } + q)}$.
\end{enumerate}

This gives us the following lemma
\begin{lemma} \label{lemma:disc-of-tower}
For the extension $T_n/T_1$, the discriminant is given by 
$$ \disc(T_n/T_1) = P_\infty^{2(q^{n-1} 
-1)}.\prod_{g_1(\alpha)=0}P_\alpha^{2(q^{n-1}-1)}.P_0^{ 2(q^{n-1} - q^{n- \lfloor \frac{n}{2} \rfloor } - q^{\lfloor \frac{n}{2} \rfloor } + q)}. $$
\end{lemma}

\section{Towards a global integral basis}

We will now look at two possible choices for the set $S$ of places of \pf. Note that \of and \oe depend on the choice of $S$, so that changing $S$ implicitly changes them also.

\subsection{$S = \pf \setminus P_\infty$}
Our first choice of $S$ is the set of all places except $P_\infty$.
\begin{lemma} 
$$ \otone = \fq[x_1]. $$
\end{lemma}
This is just saying that the only functions on the projective line with poles only at $P_\infty$ are 
the polynomials. Since $\fq[x_1]$ is a p.i.d, we know that there will be a 
global integral basis in $T_n$, over $\fq[x_1]$. Also, from Proposition~\ref{proposition:global=local-everywhere} and Theorem~\ref{theorem:setdisc=fielddisc}, we know that the discriminant of such a basis would have to be the field 
discriminant, at all places other than infinity.

We now examine more closely the set of functions ${\frak S}$ that is the 
``product'' of the sets $\{1,x_i,\ldots,x_i^{q-1}\}, \; 2 \leq i \leq n$. More specifically, 
$${\frak S} = \{\prod_{i=2}^n x_i^{e_i}\;|\;0 \leq e_i \leq q-1, \forall i \}.$$

\begin{lemma} \label{lemma:disc(S)=1}
$$  |\disc({\frak S})| = 1. $$
\end{lemma}
{\bf Proof}.  
As per \cite{Wei1}, Ch. 3.7.8, $\disc({\frak S})$ can also be written as the determinant of an $(q^{n-1} \times q^{n-1})$ matrix whose $(\underline{e},\underline{f})^{th}$ element (we use the lexicographic ordering), where
$\underline{e}=(e_2, \ldots, e_n), 0 \leq e_i \leq q-1$,  is given by 
\[
Tr_{T_n/T_1}( \prod_{i=2}^n x_i^{h_i} ) , \mbox{ where } h_i
= e_i + f_i .
\]
We can expand
\[
Tr_{T_n/T_1}( \prod_{i=2}^n x_i^{h_i} )  \ = \ 
Tr_{T_{n-1}/T_1}(  \prod_{i=2}^{n-1} x_i^{h_i} 
Tr_{T_n/T_{n-1}}( x_n^{h_n} ) ) .
\]
The minimum polynomial of $x_n$ over $T_{n-1}$ is 
$$ z^q + z = \frac{x_{n-1}^q}{x_{n-1}^{q-1} + 1}.  $$  
Let $\{ u_j \}_{1\leq j \leq q}$, (with $u_1=x_n$) denote the distinct roots in
$\overline{T_1}$ of this equation. Then 
\[
Tr_{T_n/T_{n-1}} ( x_n^{h_n}) \ = \ \sum_{j=1}^q u_j^{h_n} .
\]
Let $\nu_1,\ldots,\nu_q$ be the elementary symmetric
functions of the roots $\{ u_j \}$. Thus, $\nu_i = {\bf s}_{q,i}(u_1,u_2,\ldots u_q)$.  It follows that 
$$ \nu_1=\ldots=\nu_{q-2} = 0 \mbox{  and  } \nu_{q-1} = (-1)^{q-1} . $$ 
Let $S_r$ denote the $r$th power sum 
\[
S_r \ = \ \sum_{k=1}^q u_k^r .
\]
From Newton's identities which are stated below 
\begin{eqnarray*}
S_r & = & \nu_1 S_{r-1} - \nu_2 S_{r-2} + \ldots + (-1)^r \nu_{r-1}
S_1 + (-1)^{r+1} r \nu_r, \mbox{ for } 1 \leq r \leq q, \\ 
    & = & \nu_1 S_{r-1} - \nu_2 S_{r-2} + \ldots + (-1)^{q+1}\nu_q S_{r-q}, 
	\mbox{ for } r > q,  
\end{eqnarray*}
we obtain 
\begin{eqnarray*}
S_r & = & 0, \ \mbox{ for } \ 0 \leq r < q-1 \ \mbox{ and } \ q \leq r < 2(q-1),
\\
S_{q-1} & = & (-1)^q \nu_{q-1} (q-1) = 1, \\
S_{2(q-1)} & = & (-1)^q \nu_{q-1} S_{q-1} = -1. 
\end{eqnarray*}
Since $Tr_{T_n/T_{n-1}}(x_n^{h_n}) \ = \ S_{h_n}$, it follows that 
\begin{eqnarray*} 
Tr_{T_n/T_{n-1} }( x_n^{h_n} ) & = & 0, \mbox{ if
			$h_n$ is not divisible by $q-1$, } \\ 
 & = & 1, \mbox{ if $h_n=q-1$, }  \\
 & = & -1, \mbox{ if $h_n=2(q-1)$. }  
\end{eqnarray*} 
Thus, 
$Tr_{T_n/T_{n-1}}(x_n^{h_n})$ is nonzero only if $(q-1) | h_n$ and in this case, 
\begin{eqnarray*} 
Tr_{T_n/T_{n-1}}(x_n^{h_n}) & = & \pm 1 \cdot Tr_{T_{n-1}/T_1 } (
	\prod_{i=2}^{n-1} x_i^{h_i} ), \\
 & = & \pm 1 \cdot Tr_{T_{n-2}/T_1} ( \prod_{i=2}^{n-2} x_i^{h_i} 
	Tr_{T_{n-1}/T_{n-2}}(x_{n-1}^{h_{n-1}}  ) ).
\end{eqnarray*} 
and we can argue similarly with $n$ replaced by $n-1$ and so on to
arrive finally at 
\begin{eqnarray*} 
Tr_{T_n/T_1}( \prod_{i=2}^n x_i^{h_i} )  & = & 0, \mbox{ if
some $h_i$ is not divisible by $q-1$, } \\
 & = & 1, \mbox{ if each $h_i=q-1$. }
\end{eqnarray*} 
Thus $\disc({\frak S})$ is the determinant of a ''reverse triangular'' matrix
of the form 
\[
\left[ \begin{array}{ccccc} 
	0 & 0 & \ldots  & 0 & 1 \\
	0 & 0 &  \ldots & 1 & * \\
	\ldots & \ldots & \ldots & \ldots & \ldots \\
	0 & 1 & \ldots & * & * \\
	1 & * & * & * & *
	\end{array} \right]. 
\]
whose reverse-diagonal terms all equal $1$, and the entries to the
left and above the reverse-diagonal are all zero. The statement of the lemma follows immediately. \hfill $\Box$ 

Now consider the set ${\frak S}_{g_1}$ obtained by multiplying each element in 
${\frak S}$ by $g_1$, except for the element $'1'$, \ie,  

$$ {\frak S}_{g_1} = \{1\}\cup\{\prod_{i=2}^n x_i^{e_i}\;|\;0 \leq e_i \leq q-1,\mbox{ atleast one } e_i \neq 0 \}. $$

\begin{corollary} \label{corollary:disc-of-set-g1(S)}
$$  \disc( {\frak S}_{g_1} ) = g_1^{2(q^{n-1}-1)}. $$
\end{corollary}
{\bf Proof}. Observe that $g_1$ is fixed by the $\sigma$, and apart from $'1'$ 
there are $(q^{n-1}-1)$ elements in ${\frak S}$. 
Now we get the result using the definition of the discriminant and from Lemma~\ref{lemma:disc(S)=1}. \hfill $\Box$

\begin{lemma} \label{lemma:S-is-basis-except-at-zero}
The set ${\frak S}_{g_1}$ forms a local integral basis at all places other 
than $P_\infty$ and $P_0$.
\end{lemma}
{\bf Proof}. Follows from Theorem~\ref{theorem:setdisc=fielddisc}, Lemma~\ref{lemma:disc-of-tower}, and Corollary~\ref{corollary:disc-of-set-g1(S)}. \hfill $\Box$

\subsection{$S = \pf \setminus \{P_\infty,P_0\}$}

Lemma~\ref{lemma:S-is-basis-except-at-zero} suggests that it might be easier to work with the ring of regular functions for a new choice of $S$ that excludes even the place $P_0$. We have the following lemma that describes all functions in \of, for this choice of $S$.

\begin{lemma}
$$ \of = \sq[x_1]\left[\frac{1}{x_1}\right].$$
\end{lemma}

Thus, functions in \of are polynomials in $x_1$, divided by arbitrary powers of $x_1$.

\begin{lemma} \label{lemma:g1(S)-is-basis-except-at-zero}
The set ${\frak S}_{g_1}$ forms a basis for \oe over \of, for $S = \pf \setminus \{P_\infty, P_0\}$.
\end{lemma}
{\bf Proof}. Follows from Theorem ~\ref{theorem:setdisc=fielddisc}, Lemma~\ref{lemma:disc-of-tower} and Corollary~\ref{corollary:disc-of-set-g1(S)}.\hfill $\Box$

We are now in a position to give the main theorem of this section.
\begin{theorem}
Every function in $T_n$ that has no poles at any finite place can be written as a linear combination of functions in the set ${\frak S}_{g_1}$, with coefficients that are of the form 
$p(x_1)/x_1^i$, where $p(x_1)$ is a polynomial in $x_1$. 
\end{theorem}
{\bf Proof}. Clearly, the ring of regular functions for the new choice of $S$ contains the ring of regular functions for the earlier choice of $S$. Now the result follows from Lemma~\ref{lemma:g1(S)-is-basis-except-at-zero}. \hfill $\Box$



\begin{thebibliography}{99}

\bibitem{AleDeoVijSti1}I. Aleshkinov, V. Deolalikar, P. Vijay Kumar and H. Stichtenoth, ``Towards a basis for the space of regular functions in a tower of function fields meeting the Drinfeld-Vladut bound,'' to appear, Proceedings of the Fifth Int. Conf. on Finite Fields and Applications, Augsburg, 1999.

\bibitem{BeaPie1} R. A. Beaumont and R. S. Pierce, {\em The algebraic foundations of mathematics}, Addison-Wesley, Reading, Mass., 1963.

\bibitem{Cas1} J. W. S. Cassels, {\em Local Fields}, Cambridge University Press, Cambridge, 1986.

\bibitem{CasFro1} J. W. S. Cassels and A. Fr\"ohlich, Eds., {\em Algebraic number theory}, Academic Press, 1967.

\bibitem{Che1} C. Chevalley, {\em Introduction to the theory of algebraic 
functions of one variable}, American Mathematical Society, New York, 1951.

\bibitem{DeoEst1} V. Deolalikar, ``On splitting almost all places of degree one in extensions of function fields,'' preprint.

\bibitem{DeoEst2} V. Deolalikar, ``On splitting almost all places of degree one in extensions of function fields - II,'' preprint.

\bibitem{DeoEst3} V. Deolalikar, ``Towers of function fields meeting the Drinfeld-Vladut bound, and towers in which all places of degree one split completely,'' preprint.

\bibitem{DriVla1} V. G. Drinfeld, ``Number of points of an algebraic curve,'' Functional Analysis, {\bf 17}, 1983, pp. 53-54.

\bibitem{Elk1} N. Elkies, ``Explicit modular towers,'' Proceedings of the 
Thirty Fifth Annual Allerton Conference on Communication, Control and Computing, Urbana, IL, 1997.          

\bibitem{FreKie1} E. Freitag and R. Kiehl, {\em Etale cohomology and the Weil conjecture}, Springer-Verlag, Berlin-Heidelberg-New York, 1988.

\bibitem{FrePerSti1}G. Frey, M. Perret and H. Stichtenoth, ``On the different of abelian extensions of global fields,'' Coding Theory and Algebraic Geometry, Lecture Notes in Mathematics, {\bf 1518}, 1991, pp. 26-32.

\bibitem{Fri1} M. Fried, ``Global construction of general exceptional covers with motivation for applications to encoding,'' Contemp. Math., {\bf 168}, Amer. Math. Soc., Providence, RI, 1994. 

\bibitem{FuhSti1} R. Fuhrmann and H. Stichtenoth, ``On maximal curves,'' Journal of Number Theory, {\bf 67}, 1997, pp. 29-51.

\bibitem{GarSti1} A. Garcia and H. Stichtenoth, ``A Tower of 
Artin-Schreier extensions of function fields attaining the Drinfeld-Vladut 
bound,''  Inventiones Mathematicae, {\bf 121}, 1995, pp. 211-222.

\bibitem{GarSti2} A. Garcia and H. Stichtenoth, ``On the Asymptotic 
Behaviour of Some Towers of Function Fields over Finite Fields,'' Journal of 
Number Theory, {\bf 61}, No. 2, 1996, pp. 248-273.

\bibitem{GarSti3} A. Garcia and H. Stichtenoth, ``Asymptotically 
good towers of function fields over finite fields,'' C. R. Acad. Sci. Paris, t. {\bf 322}, S\'{e}rie I, 1996, pp. 1067-1070.

\bibitem{Gop1} V. D. Goppa, ``Codes on algebraic curves,'' Soviet Math. Doklady, {\bf 24}, No. 1, 1981, pp. 170-172.

\bibitem{Iha1} Y. Ihara, ``Some remarks on the number of rational points of algebraic curves over finite fields,'' Journal of the Faculty of Science, University of Tokyo, {\bf 28}, 1981, pp. 721-724.

\bibitem{Jan1} G. Janusz, {\em Algebraic number fields}, Academic Press, New York, 1973. 

\bibitem{Lan1} S. Lang, {\em Algebra}, Addison-Wesley, 1993.  

\bibitem{Lau1} K. Lauter, ``Ray class field constructions of curves over finite fields with many rational points,'' Algorithmic Number Theory, ed. H. Cohen, Lecture Notes in Computer Science, {\bf 1122}, Springer, 1996, pp. 187-195.

\bibitem{Man1} Y. I. Manin, ``What is the maximum number of points on a curve over ${\Bbb F}_2$?,'' Journal of the Faculty of Science, University of Tokyo, {\bf 28}, 1981, pp. 715-720.

\bibitem{NieXin1} H. Niederreiter and C. P. Xing, ``Cyclotomic function fields, Hilbert class fields, and global function fields with many rational places,'' Acta Arithmetica, {\bf 79}, No. 1, 1997, pp. 59-76.

\bibitem{ParSha1} A. N. Parshin and I. R. Shafarevich, Eds., {\em Number Theory II,} Encyclopaedia of Mathematical Sciences, {\bf 62}, Springer-Verlag, Berlin-Heidelberg-New York, 1992.

\bibitem{RucSti1} H.-G. R\"{u}ck and H. Stichtenoth, ``A characterization of Hermitian function fields over finite fields,'' J. Reine Agnew. Math., {\bf 457}, 1994, pp. 185-188.

\bibitem{Sch1} R. Schoof, ``Algebraic curves over ${\Bbb F}_2$ with many rational points,'' Journal of Number Theory, {\bf 41}, 1992, pp. 6-14.

\bibitem{Sch2} R. Schoof, {\em Algebraic curves and coding theory}, Undergraduate Texts in Math., {\bf 336}, Univ. of Trento, 1990. 

\bibitem{Ser4}J.-P. Serre, {\em Local Fields}, Springer GTM, New York, 1979.

\bibitem{Ser1}J.-P. Serre, ``Sur le nombre des points rationnels d'une courbe algebrique sur un corps fini,'' C. R. Acad. Sci. Paris S\'{e}rie I Math., {\bf 296}, 1983, pp. 397-402.

\bibitem{Ser2}J.-P. Serre, ``Nombres de points des courbes alg\'{e}briques sur  ${\Bbb F}_q$, S\'{e}m. Th\'{e}orie des Nombres 1982-1983, Exp. 22, Univ. de Bordeaux I, Talence, 1983.

\bibitem{Ser3}J.-P. Serre, {\em Rational points on curves over finite fields}, Lecture Notes, Harvard University, 1985.


\bibitem{Sti3} H. Stichtenoth, ``\"{U}ber die Automorphismengruppe eines algebraischen Funktionenkörpers von Primzahlcharakteristik. I. Eine Abschätzung der Ordnung der Automorphismengruppe,'' (German) Arch. Math. (Basel) {\bf 24}, 1973, pp. 527-544.

\bibitem{Sti2} H. Stichtenoth, ``A note on Hermitian codes over $GF(q^2)$,'' IEEE transactions on Information Theory, {\bf 34}, No. 5, September 1988, pp. 1345-1348.

\bibitem{Sti1} H. Stichtenoth, {\em Algebraic Function Fields and Codes}, 
Springer Universitext, Berlin-Heidelberg-New York, 1991.


\bibitem{Sul1}F. J. Sullivan, ``$p$-torsion in the class group of curves with too many automorphisms,'' Arch. Math., {\bf 26}, 1975, pp. 253-261.

\bibitem{TsfVlaZin1} M. A. Tsfasman, S. G. Vladut, T. Zink, ``Modular curves, Shimura curves and Goppa codes, better than the Varshamov-Gilbert bound,'' Math. Nachr., {\bf 109}, 1982, pp. 21-28.

\bibitem{Wei1} E. Weiss, {\em Algebraic number theory}, McGraw-Hill, New York, 1963.


\bibitem{Xin1}C. P. Xing, ``Multiple Kummer extension and the number of prime divisors of degree one in function fields,'' Journal of Pure and Applied Algebra, {\bf 84}, 1993, pp. 85-93.


\bibitem{XinNie2}C. P. Xing and H. Niederreiter, ``A construction of low discrepancy sequences using global function fields,'' Acta Arithmetica, {\bf 73}, 1995, pp. 87-102.


\end{thebibliography}
\end{document}